\documentclass[oneside]{amsart}
\usepackage{amssymb,cite}
\numberwithin{equation}{section}

\usepackage[a4paper,width={15cm},left=3cm,bottom=2.5cm, top=2.5cm,
includeheadfoot]{geometry}

\usepackage{latexsym}
\usepackage{textcomp}
\usepackage{amsmath}
\usepackage{amsfonts}
\usepackage{amssymb}
\usepackage{mathrsfs}
\usepackage{enumerate}

\renewcommand{\a }{\alpha }
\renewcommand{\b }{\beta }
\renewcommand{\d}{\delta }

\newcommand{\pa}{{\partial}}
\newcommand{\D }{\Delta }

\newcommand{\e }{\varepsilon }
\newcommand{\g }{\gamma}

\newcommand{\G }{\Gamma}

\newcommand{\m }{\mu }
\newcommand{\n }{\nabla }
\newcommand{\vp }{\varphi }
\renewcommand{\phi}{\varphi}

\newcommand{\s }{\sigma }

\renewcommand{\t }{\tau }

\renewcommand{\O }{\Omega }

\newcommand{\ov}{\overline}

\newcommand{\be}{\begin{equation}}
\newcommand{\ee}{\end{equation}}

\newcommand{\R}{\mathbb{R}}
\newcommand{\N}{\mathbb{N}}
\renewcommand{\S}{\mathbb{S}}

\newcommand{\de}{\partial}

\newcommand{\ti}{\widetilde}

\newcommand{\ra}{{\rangle}}
\newcommand{\la}{{\langle}}

\renewcommand{\k}{\kappa}

\newcommand{\calH }{\mathcal{H}}

\newcommand{\calC }{\mathcal{C}}

\def\dlim{\displaystyle\lim}
\def\dsup{\displaystyle\sup}

\def\dinf{\displaystyle\inf}

\newtheorem{Theorem}{Theorem}[section]

\newtheorem{Lemma}[Theorem]{Lemma}
\newtheorem{Proposition}[Theorem]{Proposition}
\newtheorem{Corollary}[Theorem]{Corollary}
\newtheorem{Remark}[Theorem]{Remark}

\def\proof{\noindent{{\bf Proof. }}}
\def\square{\vbox{
    \hrule height .4pt
    \hbox{\vrule width .4pt height 7pt \kern 7pt
       \vrule width .4pt}
    \hrule height .4pt }}

\def\square{\vbox{
    \hrule height .4pt
    \hbox{\vrule width .4pt height 7pt \kern 7pt
       \vrule width .4pt}
    \hrule height .4pt }}

\def\QED{\hfill {$\square$}\goodbreak \medskip}

\def\R{{\mathbb R}}

\def\S{{\mathbb S}}

\def\div{{\rm div}}

\newcommand{\Dsm}{(-\Delta+m^2)^s}
\newcommand{\RNp}{\R^{N+1}_+}

\newcommand{\SN}{{\mathbb S}^{N-1}}
\newcommand{\weakly}{\rightharpoonup}

\newcommand{\dive }{\mathop{\rm div}}

\begin{document}
\title[Relativistic Schr\"odinger operators with a singular
potential]{Unique continuation
  properties for relativistic Schr\"odinger operators with a singular
  potential} \author[Mouhamed Moustapha Fall \and Veronica
Felli]{Mouhamed Moustapha Fall \and Veronica Felli }
\address{\hbox{\parbox{5.7in}{\medskip\noindent
      M.M. Fall\\
      African Institute for Mathematical Sciences (A.I.M.S.) of Senegal,\\
      KM 2, Route de Joal,\\
      B.P. 1418.
      Mbour, S\'en\'egal. \\[2pt]
      {\em{E-mail address: }}{\tt mouhamed.m.fall@aims-senegal.org.}\\[5pt]
      V. Felli\\
      Universit\`a di Milano
      Bicocca,\\
      Dipartimento di Ma\-t\-ema\-ti\-ca e Applicazioni, \\
      Via Cozzi
      53, 20125 Milano, Italy. \\[2pt]
      \em{E-mail address: }{\tt veronica.felli@unimib.it.}}}}


\date{December 23, 2013}

 \begin{abstract}
   \noindent Asymptotics of solutions to relativistic fractional elliptic equations
   with Hardy type potentials is established in this paper.  As a
   consequence, unique continuation properties are obtained.
 \end{abstract}

\maketitle

\section{Introduction}\label{s:int}
Let $N>2s$ with $s\in(0,1)$ and $\Omega$ be an open subset of
$\R^N$. The purpose of the present paper is to establish unique continuation
properties for the operator
\begin{equation}\label{eq:Ham}
H:=\Dsm-\frac{a(\frac{x}{|x|})}{|x|^{2s}}-h(x),
\end{equation}
where $m\geq0$,
$a\in C^1(\S^{N-1})$, and
\begin{equation}\label{eq:ipo1}
h\in C^1(\O\setminus\{0\}),
\quad |h(x)|+|x\cdot\nabla h(x)|\leq C_h|x|^{-2s+\chi} \text{ as } |x|\rightarrow 0,
\end{equation}
for some $C_h>0$ and $\chi\in(0,1)$. Answers to the problem of unique
continuation will be derived from a precise description of the
asymptotic behavior of solutions to $Hu=0$ near $0$.

From the mathematical point of view, a reason of interest in
potentials of the type $a(x/|x|)|x|^{-2s}$ relies in their criticality
with respect to the differential operator $\Dsm$; indeed,
they have the same homogeneity as the $s$-laplacian $(-\Delta)^s$,
hence they cannot be regarded as a lower order perturbation term.
The physical  interest  in the study of properties of the Hamiltonian
in \eqref{eq:Ham} is manifest in the case $s=1/2$; indeed, if $s=1/2$
and $a\equiv Ze^2$ is constant, then the Hamiltonian \eqref{eq:Ham} describes a
spin zero relativistic particle of charge $e$ and mass $m$ in the Coulomb field of an infinitely heavy nucleus
of charge $Z$, see e.g. \cite{Herbst,lieb}.

 Before going further, let us fix our notion of solutions to $Hu=0$ in an open set $\O$.
For every $\varphi\in C^\infty_c(\R^N)$ and $s\in(0,1)$, the
relativistic Schr\"odinger operator with mass $m\geq0$  is defined as
\begin{equation}\label{eq:2}
  \Dsm \varphi(x)={c_{N,s}}m^{\frac{N+2s}{2} }
P.V.\int_{\R^{N}}\frac{ \vp(x)-\vp(y)}{|x-y|^{\frac{N+2s}{2}} } K_{
\frac{N+2s}{2}}(m|x-y|)\,dy+m ^{2s} \vp(x),
\end{equation}
for every $x\in\R^N$,
where $P.V.$ indicates  that the integral is meant in the principal
value sense
and
$$
c_{N,s}=  
{2^{-(N+2s)/2+1}}\pi^{-\frac
  N2}2^{2s}\frac{s(1-s)}{\Gamma(2-s)},
$$
see Remark \ref{rem:Dsm}.
Here  $K_{\nu}$  denotes the modified Bessel function of the second
kind with order $\nu$, see appendices B and C in sections
\ref{s:ext-th} and \ref{s:Shro}.
 The Dirichlet form associated to
$\Dsm$ on $C^\infty_c(\R^N)$ is given by
\begin{align}\label{eq:3}
  (u,v)_{
    H^{s}_m(\R^N)}:&=\int_{\R^N}(|\xi|^2+m^2)^s\widehat{u}(\xi)\ov{\widehat{v}(\xi)} d\xi\\
    &\nonumber =\frac{c_{N,s}}{2}m^{\frac{N+2s}{2} }
\int_{\R^{2N}}\frac{(u(x)-u(y))(v(x)-v(y))}{|x-y|^{\frac{N+2s}{2}} } K_{ \frac{N+2s}{2}}(m|x-y|)\,dx\,dy\\
\nonumber & \hspace{1cm} +m ^{2s} \int_{\R^N}u(x)v(x)dx,
\end{align}
where $\widehat{u}$ denotes the unitary Fourier transform of $u$.
We  define $ H^{s}_m(\R^N)$ as the completion of
$C^\infty_c(\R^N)$ with respect to the norm induced by the
scalar product \eqref{eq:3}. If $m>0$, $ H^{s}_m(\R^N)$ is nothing
but the standard $ H^{s}(\R^N)$; then, we will write  $H^{s}(\R^N)$
without the subscript ``$m$''.

By a weak solution to $Hu=0$ in $\O$, we mean a function
$u\in  H^{s}_m(\R^N)$ such that
\begin{equation}\label{eq:1}
(u,\varphi)_{ H^{s}_m(\R^N)}
=\int_\Omega \bigg(\frac{a(x/|x|)}{|x|^{2s}}u(x)+h(x)u(x) \bigg)\varphi(x)\,dx,
\text{ for all }\varphi\in C^\infty_c(\Omega).
\end{equation}
We notice that the right hand side of \eqref{eq:1} is well defined in
view of  the following Hardy type  inequality due to Herbst in \cite{Herbst}  (see also 
\cite{yafaev}):
\begin{equation}\label{eq:frac_hardy}
  \Lambda_{N,s}\int_{\R^N}\frac{u^2(x)}{|x|^{2s}}\,dx\leq \int_{\R^N}
  |\xi|^{2s}|\widehat u(\xi)|^2\,d\xi \leq \|u\|_{ H^{s}_m(\R^N)}^2
  ,\quad\text{for all }u\in    H^{s}_m(\R^N),
\end{equation}
where
$$
\Lambda_{N,s}:=2^{2s}\frac{\Gamma^2\big(\frac{N+2s}{4}\big)}{\Gamma^2\big(\frac{N-2s}{4}\big)}.
$$
 A first aim of this paper  is to give a precise description of the  behavior near $0$ of
solutions to the equation $Hu=0$, from which several unique continuation properties can be derived.
The rate and the shape of $u$ can be described in terms of
the eigenvalues and the eigenfunctions of
 the following eigenvalue problem
\begin{align}\label{eq:4}
  \begin{cases}
    -\dive\nolimits_{{\mathbb S}^{N}}(\theta_1^{1-2s}\nabla_{{\mathbb
        S}^{N}}\psi)=\mu\,
    \theta_1^{1-2s}\psi, &\text{in }{\mathbb S}^{N}_+,\\[5pt]
-\lim_{\theta_1\to 0^+} \theta_1^{1-2s}\nabla_{{\mathbb
    S}^{N}}\psi\cdot {\mathbf
  e}_1=\kappa_s a(\theta')\psi,&\text{on }\partial {\mathbb S}^{N}_+,
  \end{cases}
\end{align}
where
\begin{equation}\label{eq:kappa_s}
\kappa_s=\frac{\Gamma(1-s)}{2^{2s-1}\Gamma(s)}
\end{equation}
and
\begin{align*}
  &{\mathbb S}^{N}_+=\{(\theta_1,\theta_2,\dots, \theta_{N+1})\in
  {\mathbb S}^{N}:\theta_1>0\}=\left\{\tfrac{z}{|z|}:z\in \R^{N+1},\ z\cdot {\mathbf e}_1>0\right\},
\end{align*}
with ${\mathbf e}_1=(1,0,\dots,0)$; we refer to section
\ref{sec:separ-vari-extens} for a variational formulation of \eqref{eq:4}. From classical
spectral theory (see section \ref{sec:separ-vari-extens} for the
details), if 
\begin{equation}\label{eq:5}
\mu_1(a):=\inf_{\psi\in H^1({\mathbb
  S}^{N}_+;\theta_1^{1-2s}) \setminus\{0\}}\frac{\int_{{\mathbb
    S}^{N}_+}\theta_1^{1-2s}|\nabla\psi(\theta)|^2 dS-\kappa_s
\int_{{\mathbb S}^{N-1}}a(\theta')\psi^2(0,\theta')\,d S'}{ \int_{{\mathbb
    S}^{N}_+}\theta_1^{1-2s}\psi^2(\theta)\,dS}>-\infty,
\end{equation}
then 
problem \eqref{eq:4} admits a diverging sequence of real eigenvalues
with finite multiplicity
\begin{equation*}
  \mu_1(a)\leq\mu_2(a)\leq\cdots\leq\mu_k(a)\leq\cdots,
\end{equation*}
the first one of which coincides with the infimum in \eqref{eq:5},
which is actually attained.
Throughout the present paper, we will always assume that
\begin{equation}\label{firsteig_strict_in}
\mu_1(a)>- \bigg(\frac{N-2s}2\bigg)^{\!\!2}.
\end{equation}
Our first result is the following asymptotics of solutions at the
singularity, which generalizes to the case $m>0$ an analogous result
obtained by the authors in \cite{FF} for $m=0$.
\begin{Theorem} \label{t:asym-frac}
Let $u\in  H^{s}_m(\R^N)$ be a nontrivial weak
solution to
\[
\Dsm
u(x)-\frac{a(\frac{x}{|x|})}{|x|^{2s}}u(x)-h(x)u(x)=0
\] in an open set  $\Omega\subset
\R^N$ containing the origin, with $s\in (0,1)$, $N>2s$, $m\geq 0$,
$h$ satisfying assumption \eqref{eq:ipo1}, and $a\in C^1(\SN)$. Then
there exists an eigenvalue $\mu_{k_0}(a)$ of \eqref{eq:4} and an
eigenfunction $\psi$ associated to $\mu_{k_0}(a)$ such that
\[
\t^{-\frac{2s-N}{2}-\sqrt{ \left(\frac{2s-N}{2}  \right)^2
    +\mu_{k_0}(a)  }}u(\tau x)\to
|x|^{-\frac{N-2s}{2}+\sqrt{ \left(\frac{2s-N}{2}  \right)^2
    +\mu_{k_0}(a)  }}
\psi\big(0,\tfrac{x}{|x|}\big)\quad\textrm{as } \t\to0^+,
\]
in  $C^{1,\alpha}_{\rm
    loc}(B_1'\setminus\{0\})$ for some $\a\in(0,1)$,
where $B_1':=\{x\in \R^N:|x|<1\}$,
and, in particular,
\[
\t^{-\frac{2s-N}{2}-\sqrt{ \left(\frac{2s-N}{2}  \right)^2 +\mu_{k_0}(a)  }}u\left(\tau\theta'\right)\to \psi\left(0,\theta'\right)\quad\textrm{in } C^{1,\a}(\S^{N-1})\quad\textrm{as } \t\to0^+,
\]
where $\S^{N-1}=\de \S^N_+$.
\end{Theorem}
The proof of Theorem \ref{t:asym-frac} is based on an Almgren type
monotonicity formula (see \cite{almgren,GL}) for a
Caffarelli-Silvestre type extended problem. Indeed, for every $u\in H^s(\R^N)$ there exists a unique
$w=\mathcal H(u)\in H^1(\R^{N+1}_+;t^{1-2s})$ weakly solving
\[
\begin{cases}
  -\dive(t^{1-2s}\nabla w)+m^2 t^{1-2s}w=0,&\text{in }\R^{N+1}_+,\\
w=u,&\text{on }\partial \R^{N+1}_+=\{0\}\times \R^N,
\end{cases}
\]
where $\R^{N+1}_+
=\{z=(t,x):t\in(0,+\infty),\ x\in\R^N\}$ and  $H^1(\R^{N+1}_+;t^{1-2s})$ is defined as the  completion of
$C^\infty_c(\overline{\R^{N+1}_+})$ with respect to the norm
$$
\|w\|_{H^1(\R^{N+1}_+;t^{1-2s})}=\bigg(
\int_{\R^{N+1}_+ }t^{1-2s}\Big(|\nabla w(t,x)|^2+w^2(t,x)\Big)
dt\,dx
\bigg)^{\!\!1/2}.
$$
Furthermore,
$$
-\lim_{t\to 0^+}t^{1-2s}\frac{\partial w}{\partial
  t}(x)=\kappa_s(-\Delta+m^2)^su(x), \quad\text{in }H^{-s}(\R^N),
$$
in a weak sense, see Theorem \ref{th:ext} in Appendix B. Therefore $u\in
  H^s(\R^N)$ weakly solves $H(u)=0$ in $\Omega$ in the
  sense of \eqref{eq:1} if and only if its extension $w=\mathcal H(u)$
  satisfies
\begin{equation}\label{eq:extended}
\begin{cases}
    -\dive(t^{1-2s}\nabla  w(t,x)) +m^2 t^{1-2s}w=0,&\text{in }\R^{N+1}_+,\\
w(0,x)=u(x),&\text{in }\R^{N},\\
-\lim_{t\to 0^+}t^{1-2s}\frac{\partial w}{\partial
  t}(t,x)=\kappa_s
\Big(\frac{a(x/|x|)}{|x|^{2s}}w+hw\Big), &\text{in }\Omega,
\end{cases}
\end{equation}
in a weak sense. The asymptotics provided in Theorem \ref{t:asym-frac}
follows from combining an Almgren type monotonicity formula for
problem \eqref{eq:extended} with a
blow-up analysis; see \cite{FFT,FFT2,FFT3} for the combination of such
methods to prove not only unique continuation but also the precise
asymptotics of solutions.  We also refer to \cite{CSilv,FF} for
monotonicity formulas in fractional problems.

As a  particular case of Theorem \ref{t:asym-frac}, if $a\equiv 0$ we
obtain the following result.
\begin{Corollary} \label{t:lambda0-asym-frac}
Let $\Omega$ be an open bounded subset of $\R^N$ and $u\in  H^{s}_m(\R^N)$ be a nontrivial weak
solution to
\begin{equation}\label{eq:frac_eq_unspm}
\Dsm u(x)=h(x)u(x) ,\quad\text{in }\Omega,
\end{equation}
with $s\in(0,1)$ and  $h\in
C^1(\Omega)$. Then, for every
$x_0\in\Omega$,  there exists an eigenvalue $\mu_{k_0}=\mu_{k_0}(0)$
of problem \eqref{eq:4} with $a\equiv 0$ and an eigenfunction $\psi$
associated to $\mu_{k_0}$ such that
\begin{multline}\label{eq:limint0lambda}
  \t^{-\frac{2s-N}{2}-\sqrt{ \left(\frac{2s-N}{2} \right)^2
      +\mu_{k_0}  }}u(x_0+\tau (x-x_0))\\
  \to |x-x_0|^{-\frac{N-2s}{2}+\sqrt{ \left(\frac{2s-N}{2} \right)^2
      +\mu_{k_0} }}
  \psi\big(0,\tfrac{x-x_0}{|x-x_0|}\big)\quad\textrm{as } \t\to0^+,
\end{multline}
in  $C^{1,\alpha}(\{x\in\R^N:x-x_0\in B_1'\})$.
\end{Corollary}

\noindent A relevant application of the asymptotic analysis contained in
Theorem \ref{t:asym-frac} and Corollary \ref{t:lambda0-asym-frac} is
the validity of some unique continuation principles.  A direct consequence
of Theorem~\ref{t:asym-frac} is the following {\em{strong unique
    continuation property}}, which extends to the case $m> 0$ an
analogous result obtained for $m=0$ in \cite{FF}.

\begin{Theorem}\label{t:sun}
  Suppose that all the assumptions of Theorem \ref{t:asym-frac} hold true.
Let $u\in  H^{s}_m(\R^N)$ be a weak
solution to
\[
\Dsm
u(x)-\frac{a(\frac{x}{|x|})}{|x|^{2s}}u(x)-h(x)u(x)=0
\]
 in an open set
$\Omega\subset \R^N$ containing the origin. If
$u(x)=o(|x|^n)=o(1)|x|^n$ as $|x|\to 0$ for all $n\in
  \N$, then $u\equiv 0$ in $\Omega$.
\end{Theorem}
We mention that recently some strong unique continuation properties for
fractional laplacian have been proved by several authors, see
\cite{FF,FSV,Seo,Ru, Yang}.
\noindent Corollary \ref{t:lambda0-asym-frac} allows also to prove the
following unique continuation principle from sets of positive
measures, which implies, as an interesting application,
 that the nodal
sets of eigenfunctions for $\Dsm$ have zero Lebesgue measure.
\begin{Theorem}\label{t:unspm}
Suppose that $u$ is as in Corollary \ref{t:lambda0-asym-frac}.
If $u\equiv 0$ on a set $E\subset \Omega$ of positive measure, then $u\equiv 0$ in $\Omega$.
\end{Theorem}

A direct application of  Theorem \ref{t:unspm} can be found in \cite{FSV}, where the authors  proved the case $N=1$ and $m=0$.

\begin{Remark}
  We point out that the results presented above still
  hold for the more general nonlinear problem
$$
\Dsm u= \frac{a(\frac{x}{|x|})}{|x|^{2s}}u(x)+h(x)u(x)+f(x,u),
$$
which was considered in \cite{FF} for $m=0$. Assuming that
\begin{align}
\begin{array}{l}
f\in C^1(\Omega\times \R),
\quad t\mapsto F(x,t)\in C^{1}(\Omega\times\R), \\[5pt]
|f(x,t)t|+|f'_t(x,t)t^2|+|\nabla_x F(x,t)\cdot x|\leq C_f\,|t|^{p}
\text{ for a.e. $x\in\Omega$ and all $t\in\R$},
\end{array}
\end{align}
where $2<p\leq 2^*\!(s)=\frac{2N}{N-2s}$,
$F(x,t)=\int_0^t f(x,r)\,dr$, the asymptotics of Theorem
\ref{t:asym-frac} and the unique continuation principles of Theorems
\ref{t:sun} and \ref{t:unspm} still hold. Since the presence of the
nonlinear term introduces essentially the same difficulties already
treated in \cite{FF}, we present here the details of proofs only for
the linear problem focusing on the differences from \cite{FF} due to
the introduction of the relativistic correction.
\end{Remark}

Beside the above unique continuation properties (UCPs), several
results of independent interest will be proved in this paper. Indeed,
to prove the UCPs, we transform, in the spirit of \cite{FF}, 
problems of the type 
\begin{equation}\label{eq:6}
\Dsm
u(x)=G(x,u),\quad\text{in }\Omega,
\end{equation}
into the problem
\be\label{eq:ext-intro}
\begin{cases}
- \div(t^{1-2s}\n w)+m^2t^{1-2s} w=0,&\text{in }\RNp,\\
w=u, &\text{on }\R^N,\\
 -\lim_{t\to 0} t^{1-2s}\frac{\de w}{\de t}=\k_s\Dsm u=\k_s G(x,w), &\text{in } 
 \Omega.
\end{cases}
\ee Such extension is a generalization of the Caffarelli-Silvestre
extension \cite{CSilv} and it is a particular case of more general
extension theorems proved in Section \ref{s:ext-th}.  We actually
derive asymptotics of solutions and unique continuation for problems
of type \eqref{eq:6} as a consequence of asymptotics
and unique continuation for the corresponding extended problem
\eqref{eq:ext-intro}. 

 In sections \ref{sec:separ-vari-extens},
\ref{s:Schauder} and \ref{sec:pohoz-type-ident} we present some
preliminary results including some Hardy type inequalities, Schauder
estimates for boundary value problems related to \eqref{eq:ext-intro}
and a Pohozaev type identity. These latter preparatory results will be used in the
study of the monotonicity properties of the Almgren type frequency
function associated to the extended problem \eqref{eq:extended}; in
section \ref{sec:monot-prop} a blow-up analysis of the extended
problem will be also performed thus leading to the proof of Theorem
\ref{t:asym-frac} and, as consequences of Theorem \ref{t:asym-frac},
of Corollary \ref{t:lambda0-asym-frac} and Theorems \ref{t:sun} and
\ref{t:unspm}.  Finally, in Section \ref{s:Shro} we describe some
properties of the relativistic Schr\"{o}dinger operator $\Dsm$.

\section{Hardy type inequalities}\label{sec:separ-vari-extens}

Let us denote,
for every $R>0$,
\begin{align*}
&B_R^+=\{z=(t,x)\in\R^{N+1}_+\,:\,  |z|<R\} ,\quad
B_R':=\{x\in \R^N:|x|<R\},\\
&S_R^+=\{z=(t,x)\in\R^{N+1}_+ \,:\, |z|=R\}.
\end{align*}
For every $R>0$, we  define the space
$H^1(B_R^+;t^{1-2s})$ as the  completion of
$C^\infty(\overline{B_R^+})$ with respect to the norm
$$
\|w\|_{H^1(B_R^+;t^{1-2s})}=\bigg(
\int_{B_R^+}t^{1-2s}\Big(|\nabla w(t,x)|^2+w^2(t,x)\Big)
dt\,dx
\bigg)^{\!\!1/2}.
$$
We also define $H^{1}({\mathbb S}^{N}_+;\theta_1^{1-2s})$ as the completion of
$C^\infty(\overline{{\mathbb S}^{N}_+})$ with respect to the norm
$$
\|\psi\|_{H^1({\mathbb S}^{N}_+;\theta_1^{1-2s})}=\bigg(
\int_{{\mathbb S}^{N}_+}\theta_1^{1-2s}\big(|\nabla_{{\mathbb
        S}^{N}}\psi(\theta)|^2+\psi^2(\theta)\big)dS
\bigg)^{\!\!1/2}
$$
and
$$
L^2({\mathbb
  S}^{N}_+;\theta_1^{1-2s}):=\Big\{\psi:{\mathbb S}_+^{N}\to\R\text{
  measurable such that }{\textstyle{\int}}_{{\mathbb
    S}^{N}_+}\theta_1^{1-2s}\psi^2(\theta)\,dS<+\infty\Big\}.
$$
We recall the Sobolev trace inequality: there exists $S_{N,s}>0$ such
that, for all $w\in  \mathcal D^{1,2}(\R^{N+1}_+;t^{1-2s})$ 
$$
\bigg(\int_{\R^N}|w(0,x)|^{2N/(N-2s)}\,dx\bigg)^{\!\!(N-2s)/N}dx\leq S_{N,s}
\int_{\RNp}t^{1-2s}|\n w(t,x)|^2dt\,dx,
$$
where 
$\mathcal D^{1,2}(\R^{N+1}_+;t^{1-2s})$ 
is defined as the  completion of
$C^\infty_c(\overline{\R^{N+1}_+})$ with respect to the norm
$\big(
\int_{\R^{N+1}_+ }t^{1-2s}|\nabla w(t,x)|^2
dt\,dx
\big)^{1/2}$ (see e.g. \cite{FF} for details).
Using a change of variables and writing $w(z)=f(|z|)\psi(z/|z|)$, with
$f\in C^\infty_c(0,\infty) $, we can easily prove that there exists
a well defined continuous trace operator
\[
 H^1({\mathbb S}^{N}_+;\theta_1^{1-2s})\to L^{2N/(N-2s)}(\partial{\mathbb
    S}^{N}_+)=L^{2N/(N-2s)}({\mathbb S}^{N-1}),
\]
 so that, for some $C_{N,s}>0$,
\be\label{eq:Sob-SN}
\| \psi(0,\cdot)     \|_{  L^{2N/(N-2s)}({\mathbb S}^{N-1})}^2\leq
      C_{N,s} \left(\int_{{\mathbb S}^{N}_+}\theta_1^{1-2s}|\nabla
\psi(\theta)|^2 \,dS+ \int_{{\mathbb S}^{N}_+}\theta_1^{1-2s}
\psi^2(\theta) \,dS   \right)
 \ee
for all $\psi\in H^1({\mathbb S}^{N}_+;\theta_1^{1-2s})$.

In order to construct an
orthonormal basis of $L^2({\mathbb S}^{N}_+;\theta_1^{1-2s})$ for
expanding solutions to $Hu=0$ in Fourier series, we are
naturally lead to consider the eigenvalue problem \eqref{eq:4}, which
admits the following variational formulation: we say that $\mu\in\R$
is an eigenvalue of problem \eqref{eq:4} if there exists $\psi\in H^1({\mathbb
  S}^{N}_+;\theta_1^{1-2s})\setminus\{0\}$ (called eigenfunction) such
that
$$
Q(\psi,\vartheta)=\mu\int_{{\mathbb
    S}^{N}_+}\theta_1^{1-2s}\psi(\theta)\vartheta(\theta)\,dS,\quad\text{for
all }\vartheta \in H^1({\mathbb S}^{N}_+;\theta_1^{1-2s}),
$$
where
\begin{align*}
&Q:H^1({\mathbb S}^{N}_+;\theta_1^{1-2s})\times H^1({\mathbb
  S}^{N}_+;\theta_1^{1-2s})\to\R,\\
&Q(\psi,\vartheta)=\int_{{\mathbb S}^{N}_+}\theta_1^{1-2s}\nabla
\psi(\theta)\cdot\nabla\vartheta(\theta)\,dS-
    \kappa_s \int_{{\mathbb
      S}^{N-1}}a(\theta')\psi(0,\theta') \vartheta(0,\theta')\,dS' .
\end{align*}
  If $a\in L^{N/(2s)}({\mathbb S}^{N-1})$ and \eqref{eq:5} holds, then we can prove that  the bilinear form $Q$ is continuous and
weakly coercive on $H^1({\mathbb S}^{N}_+;\theta_1^{1-2s})$.
 Moreover,
since the weight $\theta_1^{1-2s}$ belongs to the second Muckenhoupt
class, the embedding
\[
H^1({\mathbb
  S}^{N}_+;\theta_1^{1-2s})\hookrightarrow \hookrightarrow
L^2({\mathbb S}^{N}_+;\theta_1^{1-2s})
\]
is compact. From classical spectral theory, problem \eqref{eq:4}
admits a diverging sequence of real eigenvalues with finite
multiplicity $\mu_1(a)\leq\mu_2(a)\leq\cdots\leq\mu_k(a)\leq\cdots$
the first of which coincides with the infimum in \eqref{eq:5} and then
admits the variational characterization
\begin{equation}\label{firsteig}
  \mu_1(a)=\min_{\psi\in H^1({\mathbb
  S}^{N}_+;\theta_1^{1-2s}) \setminus\{0\}}\frac{Q(\psi,\psi)}{\int_{{\mathbb
    S}^{N}_+}\theta_1^{1-2s}\psi^2(\theta)\,dS}.
\end{equation}
We assume that \eqref{firsteig_strict_in} holds.
To each $k\geq 1$, we
associate an $L^2({\mathbb
  S}^{N}_+;\theta_1^{1-2s})$-normalized
eigenfunction $\psi_k\in H^1({\mathbb
  S}^{N}_+;\theta_1^{1-2s})$, $\psi_k\not\equiv 0$ corresponding to
 the $k$-th eigenvalue $\mu_{k}(a)$,
i.e. satisfying
\begin{equation}\label{angular}
Q(\psi_k,\vartheta)=\mu_k(a)\int_{{\mathbb
    S}^{N}_+}\theta_1^{1-2s}\psi_k(\theta)\vartheta(\theta)\,dS,\quad\text{for
all }\vartheta \in H^1({\mathbb S}^{N}_+;\theta_1^{1-2s}).
\end{equation}
In the enumeration
$\mu_1(a)\leq\mu_2(a)\leq\cdots\leq\mu_k(a)\leq
\cdots$, we repeat each eigenvalue as many times as its multiplicity;
thus exactly one eigenfunction $\psi_k$ corresponds to each index
$k\in{\mathbb{N}}$, $k\geq1$. We can choose the functions $\psi_k$ in
such a way that they form an orthonormal basis of $L^2({\mathbb
  S}^{N}_+;\theta_1^{1-2s}) $.

The following results will be useful to prove Hard-type inequalities
for the potential $a(x/|x|)|x|^{-2s}$ with $a$ belonging to some $L^p$
space; indeed, the Hardy inequality for this
potential involves only $\mu_1(a)$ whose corresponding
eigenfunction is simple.

\begin{Lemma}\label{lem:a_n}
  If $a\in L^{N/2s}({\mathbb S}^{N-1})$ and $a$ satisfies
  \eqref{eq:5}, then $\mu_1(a)$ is attained by a positive
  minimizer. Moreover, the mapping $a\mapsto \mu_1(a)$ is continuous in
  $L^q(\S^{N-1})$ for every $q>N/(2s)$.
\end{Lemma}
\proof The first assertion is classical thanks to the Sobolev-trace
inequality on $\S^N_+$ \eqref{eq:Sob-SN}, so we skip the details.  Now
let $q>N/(2s)$ and $a_n\in L^q(\S^{N-1})$ such that
$a_n\to a$ in $L^q({\mathbb S}^{N-1})$ (and $a_n,a$ satisfy \eqref{eq:5}). For every $\psi \in
C^\infty(\overline{{\mathbb S}^{N}_+})$, $\psi\not\equiv0$, using H\"{o}lder inequality, we can see that
\begin{align*}
\mu_1(a_n)&\leq \frac{   \int_{{\mathbb
S}^{N}_+}\theta_1^{1-2s}|\nabla \psi(\theta)|^2 \,dS-
    \kappa_s \int_{{\mathbb
      S}^{N-1}}a_n(\theta')\psi^2(0,\theta')  \,dS'  }{ \int_{{\mathbb
S}^{N}_+}\theta_1^{1-2s} \psi^2(\theta) \,dS}\\
&\leq \frac{   \int_{{\mathbb S}^{N}_+}\theta_1^{1-2s}|\nabla
\psi(\theta)|^2 \,dS-
    \kappa_s \int_{{\mathbb
      S}^{N-1}}a(\theta')\psi^2(0,\theta')  \,dS' }{ \int_{{\mathbb
S}^{N}_+}\theta_1^{1-2s} \psi^2(\theta) \,dS}\\
&\qquad +\kappa_s\frac{  \|a_n-a\|_{L^{N/(2s)}({\mathbb
      S}^{N-1}  )}\|\psi(0,\cdot)\|^2_{  L^{2N/(N-2s)}({\mathbb
      S}^{N-1}  )   } }{   \int_{{\mathbb
S}^{N}_+}\theta_1^{1-2s} \psi^2(\theta) \,dS    }.
\end{align*}
So, choosing  $\psi$ to be a minimizer for $\mu_1(a)$, we get
$$
\mu_1(a_n)\leq \mu_1(a)+o(1),\quad\text{as }n\to\infty.
$$
 Define $\calC_\delta=\{\s\in \S^{N}_+\,:\, \textrm{dist}(\s,\de
\S^N_+)<\d\}$ for all $\d>0$. Let $\chi_\d\in C^\infty(\S^{N})$ be such
that $\chi_\d=1$ on $\calC_\delta$ and $\chi_\d=0$ on
$\S^N\setminus\calC_{2\d}$. Next, let $\psi_n$ be a positive minimizer for
$\mu_1(a_n)$ normalized so that $   \int_{{\mathbb
S}^{N}_+}\theta_1^{1-2s} \psi^2_n(\theta) \,dS=1$. Then
$$
 \begin{cases}
    -\dive\nolimits_{{\mathbb S}^{N}}(\theta_1^{1-2s}\nabla_{{\mathbb
        S}^{N}}\psi_n)=\mu_1(a_n)\,
    \theta_1^{1-2s}\psi_n, &\text{in }{\mathbb S}^{N}_+,\\[5pt]
-\lim_{\theta_1\to 0^+} \theta_1^{1-2s}\nabla_{{\mathbb
    S}^{N}}\psi_n\cdot {\mathbf
  e}_1=\kappa_s a_n(\theta')\psi_n,&\text{on }\partial {\mathbb S}^{N}_+.
  \end{cases}
$$
Multiply the above equation by $\psi_n \chi_\d^2$ and integrate by
parts  to get
\begin{align*}
  \int_{{\mathbb S}^{N}_+}\theta_1^{1-2s} &\chi_\d^2|\n \psi_n|^2(\theta)
\,dS+  2 \int_{{\mathbb S}^{N}_+}\theta_1^{1-2s} \n\psi_n\cdot
\chi_\d\psi_n\n\chi_\d(\theta) \,dS\\
&\leq  \mu_1(a)+o(1)+ \kappa_s \int_{{\mathbb
      S}^{N-1}}a_n\chi_\d^2\psi^2_n(0,\theta')  \,dS'.
      \end{align*}
Hence by H\"{o}lder's inequality
\begin{align*}
\int_{{\mathbb S}^{N}_+}\theta_1^{1-2s} & |\n
(\chi_\d\psi_n)|^2(\theta) \,dS - \int_{{\mathbb
S}^{N}_+}\theta_1^{1-2s}
|\n\chi_\d|^2 \psi_n^2(\theta) \,dS   \\
& \leq\mu_1(a)+o(1)+\k_s \|a_n\|_{L^{N/2s}(\calC_{2\d})}\|\chi_d \psi_n(0,\cdot)     \|^2_{  L^{2N/(N-2s)}({\mathbb
      S}^{N-1}  )   }
      \end{align*}
and thus
\begin{align*}
&\int_{{\mathbb S}^{N}_+}\theta_1^{1-2s}  |\n
(\chi_\d\psi_n)|^2(\theta) \,dS\\
 &\leq
       C(a,N,s,\delta)\Big(1+ \|a_n\|_{L^{q}({\mathbb
      S}^{N-1}  )}|\calC_{2\d}|^{(2sq-N)/2sq} \| \chi_\d\psi_n(0,\cdot)     \|_{  L^{2N/(N-2s)}({\mathbb
      S}^{N-1}  )   }^2\Big)
\end{align*}
for some positive constant $C(a,N,s,\delta)$ depending only on
$a,N,s,\delta$.
Therefore, provided  $\d$ is small, by the Sobolev inequality we infer
$$
\int_{\calC_\d}\theta_1^{1-2s} |\n \psi_n|^2(\theta)
\,dS=  \int_{\calC_\d}\theta_1^{1-2s} |\n (\chi_\d\psi_n)|^2(\theta)
\,dS\leq 2C(a,N,s,\d)\quad \text{for all } n\in\N.
$$
Similar arguments can be performed on geodesic balls of $\S^N_+$
with radius $\d$. By covering $\ov{\S^N_+\setminus \calC_{\d/2}}$
with such finite small balls and with a classical argument of
partition of unity, we conclude that
$$
  \int_{\S^N_+}\theta_1^{1-2s} |\n \psi_n|^2(\theta)
\,dS\leq {\rm const},\qquad |\m_1(a_n)|\leq {\rm const}.
$$
 It turns out that, up to subsequences, $ \psi_n$ converges weakly
in $ H^{1}({\mathbb S}^{N}_+;\theta_1^{1-2s})$ and strongly in
$L^2({\mathbb S}^{N}_+;\theta_1^{1-2s})$ to some nontrivial function
$\psi$, which can be easily proved to be the positive (or negative) normalized  eigenfunction associated
to $\mu_1(a)$; it then follows easily that the convergence holds for
all the sequence (not only up to subsequences) and that $\mu_1(a_n)\to
\mu_1(a)$ as $n\to \infty$.
 \QED

\begin{Lemma}\label{l:hardyboundary}
  If $a\in L^{q}({\mathbb S}^{N-1})$, with $q>N/(2s)$, then
  \begin{multline}\label{eq:10}
    \int_{B_r^+} t^{1-2s}|\nabla w|^2\,dt\,dx- \kappa_s \int_{B_r'}
    \frac{a(x/|x|)}{|x|^{2s}}w^2\,dx +
    \frac{N-2s}{2r}\int_{S_r^+}t^{1-2s}w^2\,dS\\
    \geq \left(\mu_1(a)+\bigg(\frac{N-2s}{2}\bigg)^{\!\!2}\right)
    \int_{B_r^+}t^{1-2s}\frac{w^2}{|z|^2}\,dz
     \end{multline}
    for all $r>0$ and $w\in H^1(B_r^+;t^{1-2s})$.
  \end{Lemma}
\proof
  By scaling, it is enough to prove the inequality for $r=1$.  Let
  $w\in C^\infty(\overline{B_1^+})$. We have that
\begin{align}\label{eq:4du}
  \int_{B_1^+} &t^{1-2s}|\nabla w|^2\,dt\,dx- \kappa_s \int_{B_1'}
    \frac{a(x/|x|)}{|x|^{2s}}w^2\,dx +
    \frac{N-2s}{2}\int_{S_1^+}t^{1-2s}w^2\,dS\\
\notag&=
\int_{B_1^+}t^{1-2s}\bigg(\nabla w(z)\cdot\frac{z}{|z|}\bigg)^{\!\!2}dz
+    \bigg(\frac{N-2s}{2}\bigg)\int_{S_1^+}t^{1-2s}w^2dS
\\
\notag&\quad\ +\int_0^1
  \frac{\rho^{N+1-2s}}{\rho^2}\bigg(\int_{{\mathbb
      S}^{N}_+}\theta_1^{1-2s}|\nabla_{{\mathbb S}^{N}}
  w(\rho\theta)|^2
\,dS
-\kappa_s\int_{{\mathbb S}^{N-1}}a(\theta')w^2(\rho\theta ')
\bigg)d\rho.
\end{align}
From \cite[Lemma 2.4]{FF} we have that
\begin{equation}\label{eq:2du}
\int_{B_1^+}t^{1-2s}\bigg(\nabla w(z)\cdot\frac{z}{|z|}\bigg)^{\!\!2}dz
+    \bigg(\frac{N-2s}{2}\bigg)\int_{S_1^+}t^{1-2s}w^2dS
\geq \bigg(\frac{N-2s}{2}\bigg)^{\!\!2}\int_{B_1^+}t^{1-2s}\frac{w^2}{|z|^2}\,dz,
\end{equation}
whereas, from \eqref{firsteig} it follows that
\begin{align}
  \label{eq:3du}
 &\int_0^1
  \frac{\rho^{N+1-2s}}{\rho^2}\bigg(\int_{{\mathbb
      S}^{N}_+}\theta_1^{1-2s}|\nabla_{{\mathbb S}^{N}}
  w(\rho\theta)|^2
\,dS
-\kappa_s\int_{{\mathbb S}^{N-1}}a(\theta')w^2(\rho\theta ')
\bigg)d\rho \\
\notag&\quad \geq \mu_1(a) \int_0^1
  \frac{\rho^{N+1-2s}}{\rho^2}\bigg(\int_{{\mathbb
      S}^{N}_+}\theta_1^{1-2s}w^2(\rho\theta)dS\bigg)\,d\rho=\mu_1(a) \int_{B_1^+}t^{1-2s}\frac{w^2}{|z|^2}\,dz.
\end{align}
The conclusion follows from \eqref{eq:4du}, \eqref{eq:2du}, and
\eqref{eq:3du} and density
of $C^\infty(\overline{B_1^+})$ in $H^1(B_1^+;t^{1-2s})$.
\QED

\begin{Corollary}\label{c:hardyboundary_con}
  If $a\in L^{q}({\mathbb S}^{N-1})$, with $q>N/(2s)$, satisfies \eqref{firsteig_strict_in}, then there exists
  $C_{a,N,s}>0$ such that
  \begin{multline*}
    \int_{B_r^+} t^{1-2s}|\nabla w|^2\,dt\,dx- \kappa_s \int_{B_r'}
    \frac{a(x/|x|)}{|x|^{2s}}w^2\,dx +
    \frac{N-2s}{2r}\int_{S_r^+}t^{1-2s}w^2\,dS\\
    \geq C_{a,N,s} \bigg(
    \int_{B_r^+} t^{1-2s}|\nabla w|^2\,dt\,dx+
    \frac{N-2s}{2r}\int_{S_r^+}t^{1-2s}w^2\,dS\bigg)
     \end{multline*}
    for all $r>0$ and $w\in H^1(B_r^+;t^{1-2s})$.
  \end{Corollary}
  \proof
  By scaling, it is enough to prove the inequality for $r=1$. We argue by contradiction and
  assume that, for every $\e>0$ there exists $w_\e\in  H^1(B_1^+;t^{1-2s})$
such that
 \begin{multline*}
    \int_{B_1^+} t^{1-2s}|\nabla w_\e|^2\,dt\,dx- \kappa_s \int_{B_1'}
    \frac{a(x/|x|)}{|x|^{2s}}w_\e^2\,dx +
    \frac{N-2s}{2}\int_{S_1^+}t^{1-2s}w_\e^2\,dS\\
    < \e \bigg(
    \int_{B_1^+} t^{1-2s}|\nabla w_\e|^2\,dt\,dx+
    \frac{N-2s}{2}\int_{S_1^+}t^{1-2s}w_\e^2\,dS\bigg)
     \end{multline*}
 i.e.
 \begin{equation*}
    \int_{B_1^+} t^{1-2s}|\nabla w_\e|^2\,dt\,dx +
    \frac{N-2s}{2}\int_{S_1^+}t^{1-2s}w_\e^2\,dS
- \kappa_s \int_{B_1'}
    \frac{(1-\e)^{-1}a(x/|x|)}{|x|^{2s}}w_\e^2\,dx <0.
  \end{equation*}
  From Lemma \ref{l:hardyboundary} it follows that
\[
\bigg(\mu_1\Big(\frac{a}{1-\e}\Big)+\Big(\frac{N-2s}{2}\Big)^{\!2}\bigg)
    \int_{B_1^+}t^{1-2s}\frac{w_\e^2}{|z|^2}\,dz<0
\]
and hence
\[
\mu_1\Big(\frac{a}{1-\e}\Big)+\Big(\frac{N-2s}{2}\Big)^{\!2}<0
\]
On the other hand, from Lemma \ref{lem:a_n},  letting $\e\rightarrow 0$,
  we obtain $\mu_1(a)\leq -\big(\frac{N-2s}{2}\big)^{2}$, thus
  contradicting assumption \eqref{firsteig_strict_in}.
\QED

\noindent The following corollary follows from  Proposition \ref{prop:trace} and Corollary \ref{c:hardyboundary_con}.

\begin{Corollary}\label{c:hardfou}
  If $a\in L^{q}({\mathbb S}^{N-1})$,  with $q>N/(2s)$, satisfies
  \eqref{firsteig_strict_in} and $C_{a,N,s}>0$ is as in Corollary
  \ref{c:hardyboundary_con}, then
\begin{equation}\label{eq:7}
  \int_{\R^N} |\xi| ^{2s}\widehat{ u}^2\,d\xi-   \int_{\R^N}
  \frac{a(x/|x|)}{|x|^{2s}}u^2\,dx
  \geq C_{a,N,s}
  \int_{\R^N} |\xi| ^{2s}\widehat{ u}^2\,d\xi
\end{equation}
   for all   $u\in H^s_0(\R^N)$.
  \end{Corollary}

\begin{Remark}\label{rem:cont_cost}
We notice that, if $q>N/(2s)$, then the best constant  in 
inequality \eqref{eq:7} depends continuously on $a\in L^q({\mathbb
  S}^{N-1})$. Indeed, if $C_{a,N,s}$ is  the best constant  in 
\eqref{eq:7}, arguing as in \cite[Lemma 1.1]{terracini96} and
exploiting the compactness of the map $H^1({\mathbb
  S}^{N}_+;\theta_1^{1-2s})\to\R$, $\psi\mapsto \int_{{\mathbb
    S}^{N-1}}a\psi^2$ (which easily follows from \eqref{eq:Sob-SN}),
we obtain that
\begin{align*}
C_{a,N,s}&=\inf_{\mathcal
  D^{1,2}(\R^{N+1}_+;t^{1-2s})\setminus\{0\}}\frac{
\int_{\R^{N+1}_+} t^{1-2s}|\nabla w|^2\,dt\,dx- \kappa_s \int_{\R^N}
    |x|^{-2s} a(x/|x|)w^2\,dx}{
\int_{\R^{N+1}_+} t^{1-2s}|\nabla w|^2\,dt\,dx}\\
&=1-\sup_{\mathcal
  D^{1,2}(\R^{N+1}_+;t^{1-2s})\setminus\{0\}}\frac{\kappa_s \int_{\R^N}
   |x|^{-2s}a(x/|x|)w^2\,dx}{
\int_{\R^{N+1}_+} t^{1-2s}|\nabla w|^2\,dt\,dx}\\
&=1-\max_{\psi\in H^1({\mathbb
  S}^{N}_+;\theta_1^{1-2s}) \setminus\{0\}}\frac{\kappa_s
\int_{{\mathbb S}^{N-1}}a(\theta')\psi^2(0,\theta')\,d S'}{ \int_{{\mathbb
    S}^{N}_+}\theta_1^{1-2s}|\nabla\psi(\theta)|^2 dS+\big(\frac{N-2s}{2}\big)^2\int_{{\mathbb
    S}^{N}_+}\theta_1^{1-2s}\psi^2(\theta)\,dS}.
\end{align*}
From the above characterization of $C_{a,N,s}$  it is then easy to
prove that, if $a_n\to a$ in $L^q({\mathbb S}^{N-1})$, then
$C_{a_n,N,s}\to C_{a,N,s}$ as $n\to+\infty$.
\end{Remark}

\noindent
Combining Corollary \ref{c:hardyboundary_con} with \cite[Lemma
2.5]{FF}
we obtain the following estimate.
\begin{Corollary}\label{c:hardyboundary_con2}
  If $a\in L^{q}({\mathbb S}^{N-1})$, with $q>N/(2s)$, satisfies \eqref{firsteig_strict_in},  there exists
  $C'_{a,N,s}>0$ such that
  \begin{equation*}
    \int_{B_r^+} t^{1-2s}|\nabla w|^2\,dt\,dx- \kappa_s \int_{B_r'}
    \frac{a(x/|x|)}{|x|^{2s}}w^2\,dx +
    \frac{N-2s}{2r}\int_{S_r^+}t^{1-2s}w^2\,dS
    \geq C'_{a,N,s}
\int_{B_r'}\frac{
  w^2}{|x|^{2s}}\,dx
\end{equation*}
for all $r>0$ and $w\in H^1(B_r^+;t^{1-2s})$.
  \end{Corollary}

 \section{Schauder estimates for degenerate elliptic equations}\label{s:Schauder}

 As stated in Section \ref{s:int}, for $u\in H^s(\R^N)$, the nonlocal equation 
 $$
 \Dsm u=G(x,u), \quad \text{in }\Omega,
 $$
 can be reformulated as a local problem by considering its extension  in $\R^{N+1}_+$. Indeed,
  letting  $w\in H^{1}(\RNp;t^{1-2s}) $  be the unique weak solution  to the problem
$$
\begin{cases}
- \div(t^{1-2s}\n w)+m^2t^{1-2s} w=0,&\text{in }\RNp,\\
w=u, &\text{on }\R^N,
\end{cases}
$$
we have that 
$$
 -\lim_{t\to 0} t^{1-2s}\frac{\de w}{\de t}=\k_s\Dsm u,\quad \textrm{in }  \Omega,
$$
in a weak sense. This will be proved in the appendix A. This naturally
leads to the study of regularity properties of solutions to
$$
\begin{cases}
- \div(t^{1-2s}\n w)+m^2t^{1-2s} w=F(x,w),&\text{in } \Omega\times(0,T),\\
 -\lim_{t\to 0} t^{1-2s}\frac{\de w}{\de t}=G(x,w), &\text{in }\Omega,
\end{cases}
$$
which is the content of this section.

Before going on, let us state the following weighted Sobolev
inequality whose proof is essentially contained in the book of Opic
and Kufner, \cite{OK}.
\begin{Lemma}
Let $N>2s$. Then there exists a constant $S_{N,s}>0$ such that
 for every $v\in C^1_c(\R^{N+1})$ we have
\be\label{eq:wSob}
\left( \int_{\R^{N+1}_+}t^{1-2s}|v|^{\frac{2N_s}{N_s-2}}\,dt\,dx\right)^{\!\!\frac{N_s-2}{N_s}}
\leq S_{N,s}  \int_{\R^{N+1}_+}t^{1-2s}|\n v|^2 \,dt\,dx,
\ee
where $N_s= N+2-2s$.
\end{Lemma}
\proof
We have, see \cite [Section 19]{OK}, that
$$
\left( \int_{\R^{N+1}_+}t^{1-2s}|v|^{\frac{2N_s}{N_s-2}}\,dt\,dx\right)^{\!\!\frac{N_s-2}{N_s}}
\leq C_{N,s}  \left(\int_{\R^{N+1}_+}t^{1-2s}|\n v|^2 \,dt\,dx+ \int_{\R^{N+1}_+}t^{1-2s} v^2 \,dt\,dx\right).
$$
Using simple scaling argument, we obtain  \eqref{eq:wSob}.
\QED

 We will also need the following result.
 \begin{Lemma}\label{lem:estH1}
 Let $a,b\in L^p(B_1')$, for some $p>\frac{N}{2s}$ and $c,d\in L^{q}(B_1^+;t^{1-2s})$,
for some $q>\frac{N+2-2s}{2}$.
 Let $w\in H^1(B^+_1;t^{1-2s})$ be such that
$$
\begin{cases}
-\mathrm{div}(t^{1-2s}\nabla w)+ t^{1-2s} c(z)w=  t^{1-2s} d(z),
&\text{ in } B^+_1, \\
-\dlim_{t\rightarrow 0^+}t^{1-2s}\pa_tw= a(x)w+b(x),\quad  &\text{ on } B'_1.
\end{cases}
$$
Then there exits a constant $C>0$ depending  only on $N,s, \|a\|_{L^{p}(B_1')},\|c\|_{L^q(B_1^+;t^{1-2s})}$
such that
$$
\| w\|_{ H^1(B^+_1;t^{1-2s}) }\leq C\left( \| w\|_{ L^2(B^+_1;t^{1-2s}) } +\|b\|_{L^p(B_1')}+
\| d\|_{ L^q(B^+_1;t^{1-2s}) }\right).
$$
 \end{Lemma}
 \proof
 The proof is not difficult  taking into account the weighted
 Sobolev inequality \eqref{eq:wSob} together with  the  Sobolev-trace inequality:
 for every $v\in C^1_c(\R^{N+1})$
 $$
 C_{N,s}\int_{\R^N}|v(0,x)|^{\frac{2N}{N-2s}}\,dx\leq \int_{\RNp}t^{1-2s}|\n v(t,x)|^2\,dt\,dx.
 $$
 We skip the details.
  \QED

\begin{Proposition} \label{lem3.1}
Let $a,b\in L^p(B_1')$, for some $p>N/(2s)$ and $c,d\in L^{q}(B_1^+;t^{1-2s})$,
for some $q>\frac{N+2-2s}{2}$.
Let $w\in H^1(B_1^+;t^{1-2s})$ be a weak solution of
\be\label{3.1-1}
\begin{cases}
-\mathrm{div}(t^{1-2s}\nabla w)+t^{1-2s}c(z)  w\leq t^{1-2s}d(z),
&\text{ in } B_1^+, \\
-\dlim_{t\rightarrow 0^+}t^{1-2s}\pa_t w\leq a(x)w+b(x),  &\text{ on } B_1'.
\end{cases}
\ee
Then
\[
\dsup_{B_{1/2}^+}w^+\leq C\Big(\|w^+\|_{L^{2}(B_1^+;t^{1-2s})}+\|b^+\|_{L^p( B_1')}+\|d^+\|_{L^q(B_1^+;t^{1-2s})}\Big),
\]
where $w^+=\max\{0,w\}$, and $C>0$ depends only on $N,s, \|a^+\|_{L^{p}(B_1')},\|c^-\|_{L^q(B_1^+;t^{1-2s})}$.
\end{Proposition}

\proof
Let $k=\max(\|d^+\|_{L^q(B_1^+;t^{1-2s})} ,\|b^+\|_{L^p( B_1')} )$ or an arbitrary
positive small number if $\max(\|d^+\|_{L^q(B_1^+;t^{1-2s})} ,\|b^+\|_{L^p( B_1')} )=0$.
For every $L>0$, set $\overline w=w^++k$ and
\[
  \overline w_L=
\begin{cases}
\overline{w},&\quad \mbox{if } w<L,\\
k+L, &\quad \mbox{if } w\geq L.
\end{cases}
\]
Put
\[
W= \overline{w}_L^{\frac{\b}{2}}\ov{w},\qquad \phi=\eta^2(\overline w_L^\beta\overline w-k^{\beta+1})\in H^1(B_1^+;t^{1-2s}),
\]
for some $\beta\geq 0$ and some nonnegative function $\eta\in
C^1_c(B_1^+\cup B_1')$. Following \cite{JLX,TX}, testing (\ref{3.1-1})
with $\phi$, integration by parts, we have
\begin{align}\label{eq:revS}
 \int_{B_1^+}t^{1-2s }|\nabla(\eta W)|^2&\leq (1+\beta)^{\frac{\d}{\theta}}
 C\int_{B_1^+}t^{1-2s}W^2(|\nabla \eta|^2+\eta^2)\\
 \notag&\quad+2(1+\b)\int_{B_1^+}t^{1-2s}\Big(c^-+\frac{d^+}{k}\Big)(\eta W)^2 ,
\end{align}
for some positive constants $\d,\theta$ depending only on $N,s$ and $C$
depending only on $N,s,\|a^+\|_{L^p(B_1)}$.
By using H\"{o}lder inequality, we get
\begin{align}\label{eq:estcm}
\int_{B_1^+}t^{1-2s}\Big(c^-+\frac{d^+}{k}\Big)(\eta W)^2 &\leq
(\|c^-\|_{L^q(B_1^+;t^{1-2s})}+1)\|(\eta W)^2\|_{L^{\frac{q}{q-1}}(B_1^+;t^{1-2s})}\\
\notag&=: C_1\|(\eta W)^2\|_{L^{\frac{q}{q-1}}(B_1^+)}.
\end{align}
Since $1< \frac{q}{q-1}<\frac{N+2-2s}{N-2s}$, by interpolation and Young's inequalities,  we have
$$
2C_1(1+\b)\|(\eta W)^2\|_{L^{\frac{q}{q-1}}(B_1^+)}\leq \frac{1}{2S_{N,s}} \|(\eta W)^2\|_{L^{ \frac{N+2-2s}{N-2s}}(B_1^+)}+
 (1+\beta)^{\frac{\d}{\theta}} C \|(\eta W)\|_{L^{ 1}(B_1^+)}.
$$
By the weighted Sobolev inequality \eqref{eq:wSob}, we have
$$
\|(\eta W)^2\|_{L^{ \frac{N+2-2s}{N-2s}}(B_1^+)}\leq S_{N,s}\int_{B_1^+}t^{1-2s}|\n (\eta W)|^2.
$$
Using the two inequalities above in \eqref{eq:estcm}, we get
$$
2(1+\b)\int_{B_1^+}t^{1-2s}\Big(c^-+\frac{d^+}{k}\Big)
(\eta W)^2\leq  \frac{1}{2} \int_{B_1^+}t^{1-2s}|\n (\eta W)|^2+
 (1+\beta)^{\frac{\d}{\theta}} C \|\eta W\|_{L^{ 2}(B_1^+)}^2
$$
Putting this in \eqref{eq:revS}, we obtain
\[
\int_{B_1^+}t^{1-2\sigma}|\nabla(\eta W)|^2\leq
C(1+\beta)^{\frac{\delta}{\theta}}
\int_{B_1^+}t^{1-2\sigma}(\eta^2+|\nabla \eta|^2)W^2.
\]
At this point, the argument in \cite [Proposition 3.1]{TX} yields the result.
\QED

\noindent The next result is a weak Harnack inequality.

\begin{Proposition} \label{lem3.2}
Let $a,b\in L^p(B_1')$ for some $p>N/(2s)$ and $c,d\in L^{q}(B_1^+;t^{1-2s})$
for some $q>\frac{N+2-2s}{2}$.
Let $w\in H^1(B_1^+;t^{1-2s})$ be a nonnegative weak solution of
\be\label{3.4}
\begin{cases}
-\mathrm{div}(t^{1-2s}\nabla w)+c(z)t^{1-2s}w\geq  t^{1-2s}d(z),\quad &\mbox{in } B_1^+, \\
-\dlim_{t\rightarrow 0^+}t^{1-2s}\pa_t w\geq a(x)w+b(x),\quad  &\mbox{on } B_1'.
\end{cases}
\ee
Then for some $p_0>0$ and any $0<r<r'<1$ we have that
\[
\dinf_{\overline B_{r}^+}w+\|b^-\|_{L^p(B_1')}+\|d^-\|_{L^{q}(B_1^+;t^{1-2s})} \geq C\|w\|_{L^{p_0}(t^{1-2s},B_{r'}^+)},
\]
where $C>0$ depends only on $N,s,r,r', \|a^-\|_{L^{p}(B_1')},\|c^+\|_{L^{q}(B_1^+;t^{1-2s})}$.
\end{Proposition}

\proof
Set $\overline w= w+k>0$, for some positive $k$ to be determined and $v=\overline w^{-1}$.
Let $\Phi$ be any nonnegative function in $ H^1(B_1^+;t^{1-2s })$
with compact support in $B_1^+\cup B'_1$.
Multiplying both sides of the first inequality in (\ref{3.4})  by $\overline w^{-2}\Phi$ and integrating by
parts, we obtain
\[
\int_{B_1^+}t^{1-2s}\nabla v\nabla \Phi+\int_{B_1^+}t^{1-2s } \ti{c}(z)v\Phi-\int_{B'_1}\tilde a v\Phi\leq 0,
\]
where
\[
\tilde a=\frac{a^-w+b^-}{\overline w},\qquad \tilde c=\frac{c^+w+d^-}{\overline w}.
\]
If
$\max(\|b^-\|_{L^p(B_1')},\|d^-\|_{L^{q}(B_1^+;t^{1-2s}} )\neq0 $
then we choose $k=\max(\|b^-\|_{L^p(B_1')},\|d^-\|_{L^{q}(B_1^+;t^{1-2s}} )$. Otherwise,
choose an arbitrary $k>0$ which will be sent to zero.
Therefore Proposition \ref{lem3.1} (see also \cite{JLX}) implies that for any $r'\in (r,1)$
and any $p>0$
\[
\dsup_{B_r^+} v\leq C\| v\|_{L^p(B_{r'}^+;t^{1-2s})}.
\]
Following exactly the same arguments as in \cite{TX}, we get the result.
\QED

We now prove local  Schauder estimates.
\begin{Proposition}\label{p:Hold}
Let $a,b\in L^p(B'_1)$, for some $p>\frac{N}{2s}$ and $c,d\in L^{q}(B_1^+;t^{1-2s})$, for some $q>\frac{N+2-2s}{2}$.
 Let $w\in H^1(B^+_1;t^{1-2s})$ be a weak solution of
\be\label{3.1}
\begin{cases}
-\mathrm{div}(t^{1-2s}\nabla w)+ t^{1-2s}c(z)  w=  t^{1-2s}d(z),\quad &\mbox{in } B^+_1, \\
-\dlim_{t\rightarrow 0^+}t^{1-2s}\pa_tw= a(x)w+b(x),\quad  &\mbox{on
} B'_1.
\end{cases}
\ee
Then
$w\in C^{0,\a}(\ov{B_{1/2}^+})$ and in addition
$$
\|w\|_{  C^{0,\a}(\ov{B_{1/2}^+})}\leq
C\left( \|w\|_{ L^2(B^+_{1})}+\|b\|_{ L^p(B_1')}+\|d\|_{L^q(B_1^+;t^{1-2s})}  \right),
$$
with $C,\a>0$ depending only on $N,s, \|a\|_{L^{p}(B_1')},\|c\|_{L^q(B_1^+;t^{1-2s})}$.
\end{Proposition}
\proof
The proof is  a consequence of Propositions \ref{lem3.1} and \ref{lem3.2} with a standard
scaling argument for which we refer to \cite{GT}.
\QED
\begin{Remark}\label{rem:regdiv}
 Let $w\in H^1(B^+_1;t^{1-2s})$ be a weak solution of
\be\label{3.1}
\begin{cases}
-\mathrm{div}(t^{1-2s}A(z)\nabla w)+ t^{1-2s}c(z)  w=  t^{1-2s}d(z),\quad &\mbox{in } B^+_1, \\
-\dlim_{t\rightarrow 0^+}t^{1-2s}\pa_tw= a(x)w+b(x),\quad  &\mbox{on
} B'_1,
\end{cases}
\ee
with $a,b,c,d$ as in Proposition \ref{p:Hold} and the matrix $A$
satisfying
$$
C_1|\xi |^2\leq A(z)\xi\cdot\xi\leq C_2 |\xi |^2\quad \text{for all } z\in
B^+_1,\quad \xi\in\R^N,
$$
with $C_1,C_2>0$. Then the same conclusion   as in Proposition
\ref{p:Hold} holds taking into account the constants $C_1, C_2$.
\end{Remark}
%%%%%%%%%%%%%%%%%%%%%%%%%%%%%%%%%%%%%%%%%%%%%%%%%%%%%%%%%%%%%%%5
\begin{Proposition}\label{p:Holdk}
Let $a,b\in C^k(B'_1)$ and $ \n_x^i c,\n_x^i d\in L^\infty(B_1^+)$, for some $k\geq 1$ and $i=0,\dots, k$.
 Let $w\in H^1(B^+_1;t^{1-2s})$ be a weak solution of
\be\label{3.1}
\begin{cases}
-\mathrm{div}(t^{1-2s}\nabla w)+ t^{1-2s}c(z)  w=  t^{1-2s}d(z),\quad &\mbox{in } B^+_1, \\
-\dlim_{t\rightarrow 0^+}t^{1-2s}\pa_tw= a(x)w+b(x),\quad  &\mbox{on } B'_1.
\end{cases}
\ee
Then for $i=1,\dots, k$ we have that
$w\in C^{i,\a}({B_{r}^+})$, for some $r\in(0,1)$ depending only on $k$, and in addition
$$
\begin{aligned}
\sum_{i=1}^k\|\n_x^i w\|_{  C^{0,\a}(\ov{B_{r}^+})}\leq
C\bigg(  \|  w\|_{ L^2(B^+_{1};t^{1-2s})}&+\|a\|_{C^{k}(\ov{B_{r}'})}+\|b\|_{C^{k}(\ov{B_{r}'})}
+\sum_{i=1}^k\|\n_x^i c,\n_x^i d\|_{L^\infty(B_{1}^+)}  \bigg),
\end{aligned}
$$
with $C,\a>0$ depending only on $N,s, k,r, \|a\|_{L^{\infty}(B_{1/2}')},\|c\|_{L^\infty(B_{1/2}^+)}$.
\end{Proposition}
\proof
Let $h\in\R^N$ such that $|h|<\frac{1}{2}$. Then we have
$$
\begin{cases}
-\mathrm{div}(t^{1-2s}\nabla w^h)+ t^{1-2s}c(z)  w^h= - t^{1-2s}c^h(z)w+ t^{1-2s}d^h(z),\quad &\mbox{in } B^+_{1/2} \\
-\dlim_{t\rightarrow 0^+}t^{1-2s}\pa_tw^h= a(x)w^h+ a^h(x)w+b^h(x),\quad  &\mbox{on } B'_{1/2},
\end{cases}
$$
where we denote $f^h(t,x)=\frac{f(t,x+h)-f(t,x)}{h}$, for $t\geq0$.
Applying Lemma \ref{lem:estH1}, Proposition  \ref{lem3.1}  and Proposition \ref{p:Hold}  we get
\begin{align*}
 & \|w^h\|_{H^1(B^+_{1/2};t^{1-2s})} +\|w^h\|_{C^{0,\a}(\ov{ B^+_{1/4}
    })}\\
&\quad \leq C\left( \|w^h\|_{ L^2(B^+_{1/2};t^{1-2s})}
    +\|a^hw+b^h\|_{ L^\infty(B_{1/2}')}+\|c^hw+d^h\|_{L^\infty(B_{1/2}^+;t^{1-2s})}  \right)\\
  &\quad \leq C\left( \|\n w\|_{ L^2(B^+_{1/2};t^{1-2s})} +\|w\|_{C^0(\ov{B^+_{1/2}})} +\|\n_xa,\n_xb\|_{ L^\infty(B_{1/2}')}
    +\|\n_xc,\n_x d\|_{L^\infty(B_{1/2}^+;t^{1-2s})}  \right)\\
  &\quad \leq  C\left( \|  w\|_{ L^2(B^+_{1/2};t^{1-2s})}+\sum_{i=0}^1 \|\n_x^ia,\n_x^ib\|_{ L^\infty(B_{1/2}')} +\sum_{i=0}^1\|\n_x^ic,\n_x^i d\|_{L^\infty(B_{1/2}^+)
    } \right)
\end{align*}
for some positive constant $C$ depending only on $N,s, \|a\|_{L^{\infty}(B_{1/2}')},\|c\|_{L^\infty(B_{1/2}^+)}$.
Therefore using Fatou's Lemma, we obtain $W_1:=\n_x w\in H^1(B^+_{1/2};t^{1-2s})\cap C^0(B_{1/2}^+)$,
\be\label{eq:estnsw}
\|\n_x w\|_{L^{\infty}( B^+_{1/4} )}\leq\!
C\bigg(\! \|  w\|_{ L^2(B^+_{1/2};t^{1-2s})}+\sum_{i=0}^1 \|\n_x^ia,\n_x^ib\|_{ L^\infty(B_{1/2}')}
+\sum_{i=0}^1\|\n_x^ic,\n_x^i d\|_{L^\infty(B_{1/2}^+) }\!  \bigg),
\ee
and
$$
\begin{cases}
-\mathrm{div}(t^{1-2s}\nabla W_1)+ t^{1-2s}c(z)  W_1=  t^{1-2s}d_1(z),\quad &\mbox{in } B^+_{1/2} \\
-\dlim_{t\rightarrow 0^+}t^{1-2s}\pa_t W_1 =a(x)W_1+ b_1(x),\quad  &\mbox{on } B'_{1/2},
\end{cases}
$$
where $d_1(z)=-\n_xc(z)w+  \n_xd(z)$ and $b_1(x)=  \n_xa(x)w+\n_xb(x)$.
Hence by Proposition \ref{p:Hold} and \eqref{eq:estnsw}, we have
$$
\begin{aligned}
  \|W_1\|_{ C^{0,\a}(\ov{B_{1/4}^+})}&\leq
  C\left( \|W_1\|_{ L^\infty(B^+_{1/2})}+\|b_1(x)\|_{L^\infty(B_{1/2}')}+ \|d_1(z)\|_{L^\infty(B_{1/2}^+)}  \right)\\
  &\leq C\bigg( \| w\|_{ L^2(B^+_{1};t^{1-2s})}+\sum_{i=0}^1
  \|\n_x^ia,\n_x^ib\|_{
    L^\infty(B_{1/2}')}+\sum_{i=0}^1\|\n_x^ic,\n_x^i
  d\|_{L^\infty(B_{1/2}^+ )} \bigg),
\end{aligned}
$$
with $C>0$ depending only on $N,s, \|a\|_{L^{\infty}(B_{1/2}')},\|c\|_{L^\infty(\ov{B_{1/2}^+})}$.
Iterating this procedure we get the the desired estimate for $W_i=\n_x^i w$.
\QED

\section{A Pohozaev type identity}\label{sec:pohoz-type-ident}
In order to differentiate the Almgren frequency function associated to
the extended problem (see section \ref{sec:monot-prop}), we need to
derive a  Pohozaev type identity, which first requires the following
regularity result.

\begin{Lemma}\label{lem:reg}
Let $v\in H^1(B_{1}^+;t^{1-2s})$ satisfy
 \be
\label{eq:eqreg}
\begin{cases}
- \dive(t^{1-2s}\n v)+m^2t^{1-2s}v=0, &\text{in } B_1^+,\\
-\lim_{t\to 0^+}t^{1-2s}v_t=g,  &\text{on }B_1'  ,
\end{cases}
\ee
 where $g\in C^{0,\g}(B_r')$, $\g\in[0,2-2s)$ (meaning  that
      $C^{0,\g}=L^\infty$  if $\gamma=0$).
Then for every $t_0>0$ sufficiently
small there exist   positive constants $C$
and $\a\geq0$ (with $\a>0$ if $\g>0$), depending only on
$N,s,t_0,m,\g$
 such that
\be\label{eq:estnxdtwL} \|t^{1-2s} v_t\|_{ C^{0,\a}([0,t_0)\times
B_{1/8}') } \leq  C\left(\|v\|_{L^2 (B_1^+;t^{1-2s })}+  \|
g\|_{C^{\g}(B_{1/2}')}  \right). \ee
\end{Lemma}
\proof If $m=0$, this was proved in \cite{CS}. We will assume in
the following that $m>0$.
Next pick $\eta\in C^\infty_c(\ov{B_1'})$ with  $\eta=1$ on $B_{1/2}'$
and $\eta=0$ on  $\R^N\setminus\ov{B_{3/2}'}$. Then we have  that $\eta g
\in L^2(\R^N)$. By minimization arguments, there exists $W\in H^1(\RNp;
t^{1-2s})$ satisfying
$$
\begin{cases}
- \dive(t^{1-2s}\n W)+m^2t^{1-2s}W=0, &\textrm{ in } \RNp,\\
-\lim_{t\to 0^+}t^{1-2s}W_t= \eta g, &\textrm{ on } \R^N.
\end{cases}
$$
We define $w=-t^{1-2s}W_t$ and we observe that ${w}\in
L^2(\RNp;t^{-1+2s})$ and
$$
\begin{cases}
- \dive(t^{-1+2s}\n {w})+m^2t^{-1+2s} {w}=0,&\textrm{ in } \R^{N+1}_+\\
 {w}= \eta g , &\textrm{ on } \R^{N}.
\end{cases}
$$
From Remark \ref{rem:Dsm} and Proposition~\ref{eq:Poiss-cont}, it
follows that $w=\bar{P}(t,\cdot)*( \eta g)$,
where $\bar{P}$ is the Bessel kernel for the conjugate problem given
by
$$
\bar{P}(z)=C_{N,s}'\, t^{2-2s}m^{\frac{N+2-2s}{2} } |z|^{-\frac{N+2-2s}{2}}K_{\frac{N+2-2s}{2}}(m|z|),
$$
see \eqref{eq:Bess-Kernel};
we refer to Section \ref{s:PT} for asymptotics of the Bessel
function $K_{\nu}$.

\smallskip
\noindent \textbf{Claim:}   $ w\in C^{0,\g} ( \ov{\R^{N+1}_+})$ for
every $R>0$ and
\be\label{eq:Holdestw}
\|w\|_{
C^{0,\g} ( \ov{B_R^+})}\leq C_{N,s,m,R}\|\eta g\|_{C^{0,\g}(\R^N)}.
\ee
 Indeed, by a change of variables, we have that
$$
w(t,x)=C'_{N,s} \int_{\R^N}(1+|y|^2)^{-\frac{N_s}{2}}\left((tm(1+|y|^2)^{1/2}\right)^{\frac{N_s}{2}}
K_{\frac{N_s}{2}}\left((tm(1+|y|^2)^{1/2}\right)(\eta g)(x-ty)dy,
$$
where $N_s=N+2-2s$.
Let us set $f(t,|y|)=(tm(1+|y|^2)^{1/2})^{N_s/2}K_{\frac{N_s}{2}}\left((tm(1+|y|^2)^{1/2}\right)$
and $u(x)=\eta(x)g(x) $.  Letting  $x_1,x_2\in B'_R$ and $0\leq t_2<t_1<1$, we have
\begin{align}\label{eq:u1mu2}
w(t_1,x_1)-w(t_2,x_2) &=\int_{\R^N}(1+|y|^2)^{-\frac{N_s}{2}}[f(t_1,|y|)-f(t_2,|y|)]u(x_1-t_1y)dy\\
\notag&\quad +\int_{\R^N}(1+|y|^2)^{-\frac{N_s}{2}}[u(x_1-t_1y)-u(x_2-t_2y)]f(t_2,|y|)dy.
\end{align}
Using the fact that $K_{\frac{N_s}{2}}'=-\frac{N_s}{2r}K_\frac{N_s}{2}-K_{\frac{N_s}{2}-1}$, we infer that
$$
\Big|\Big (r^\frac{N_s}{2} K_\frac{N_s}{2}\Big)' \Big|=\Big|-r^\frac{N_s}{2}
K_{\frac{N_s}{2}-1 }\Big|\leq C_{N,s}m\quad  \text{for }N>2s.
$$
It follows that
$$
|f(t_1,|y|)-f(t_2,|y|)|\leq C_{N,s,m}|t_1-t_2|(1+|y|^2)^{1/2}.
$$
We recall that $\mathop{\rm supp}u\subset B'_{3/2}$ and observe that $|y|\leq \frac{1}{t_1}(3/2+2R)\leq
\frac{2}{t_1-t_2}(3/2+2R) $
provided $|x_1-t_1y|\leq\frac{3}{2}$. Therefore
\begin{align}
\label{eq:Lipf}\int_{\R^N}(1+|y|^2)^{-\frac{N_s}{2}}&|f(t_1,|y|)-f(t_2,|y|)||u(x_1-t_1y)|dy\\
&\nonumber\leq
C_{N,s,m}|t_1-t_2|\|u\|_{L^\infty(\R^N)}\int_{\{|y|\leq
\frac{1}{t_1-t_2}(3/2+2R)\}}(1+|y|^2)^{-\frac{N_s-1}{2}}dy\\
& \nonumber\leq C_{N,s,m,R}
\|u\|_{L^\infty(\R^N)}|t_1-t_2|^{2-2s}.
\end{align}
Next we have, for  $\g\in[0,2-2s)$,
\begin{align*}
  &\int_{\R^N}(1+|y|^2)^{-\frac{N_s}{2}}|u(x_1-t_1y)-u(x_2-t_2y)||f(t_2,|y|)|dy \\
  &\leq \|u\|_{C^{0,\g}(\R^N)} \|f\|_{L^\infty(\R^+\times
    \R^+)}\bigg(|t_1-t_2|^\g\!\!\int_{\R^N}
  (1+|y|^2)^{-\frac{N_s}{2}}|y|^{\g}dy+|x_1-x_2|^\g\!\!\int_{\R^N}
  (1+|y|^2)^{-\frac{N_s}{2}} dy\bigg) .
\end{align*}
Hence, for every $\g\in[0,2-2s)$,
$$
\begin{aligned}
\int_{\R^N}(1+|y|^2)^{-\frac{N_s}{2}}|u(x_1-t_1y)&-u(x_2-t_2y)||f(t_2,|y|)|dy\\
\leq C_{N,s,m}\|u\|_{C^{0,\g}(\R^N)}&
\|f\|_{L^\infty(\R^+\times \R^+)}(|t_1-t_2|^\g+ |x_1-x_2|^\g).
\end{aligned}
$$
This, together with \eqref{eq:Lipf} in \eqref{eq:u1mu2}, proves the claim.

\smallskip
\noindent We have that  $U:=v-W $  satisfies
$$
\begin{cases}
- \dive(t^{1-2s}\n U)+m^2t^{1-2s}U=0, &\textrm{ in } B_1^+,\\
t^{1-2s}U_t= 0, &\textrm{ on } B_{1/2}'
\end{cases}
$$
and  $U:=v-W \in H^1(B_1^+; t^{1-2s})$. We deduce that, for some
positive constants $C, \b$ depending only on $N,s,m$,
$$
\begin{aligned}
\| \n_x^2U\|_{ C^{0,\b}(\ov{B_{1/4}^+})}\leq C_{s,N,m} \|U\|_{
L^2(B_1^+; t^{1-2s})} \leq C\left( \|v\|_{ L^2(B_1^+; t^{1-2s})}
+\|g\|_{L^\infty(B_{1/2}')}\right)
\end{aligned}
$$
by Proposition \ref{p:Holdk}. We also observe that
$$
-t^{1-2s}\Delta_x U- ( t^{1-2s} U_t)_t+m^2t^{1-2s}U =0,\quad\text{in
}B^+_1.
$$
Then, by integration, we obtain that, for every $x\in B_{1/4}'$ and $0<t_0,t\leq
1/4$,
\be\label{eq:int-tr}
 t^{1-2s}U_t(t,x)=
t_0^{1-2s}U_t(t_0,x)-\int_t^{t_0 } \t^{1-2s}  \Delta_x U(\t,x )d\t
+m^2\int_{t}^{t_0}\t^{1-2s} U(\t,x)d\t.
\ee
Therefore
$t^{1-2s}U_t\in C^{0,\a}( \ov{B_{1/8}^+}) $ and thus
$t^{1-2s}v_t=t^{1-2s}U_t-w\in C^{0,\a}( \ov{B_{1/8}^+}) $  from
which we deduce that $\|t^{1-2s}v_t \|_{ C^{0,\a}( \ov{B_{1/8}^+})
}\leq C_{s,N,m} \left( \|v\|_{ L^2(B_1^+; t^{1-2s})}
+\|g\|_{C^\g(B_{1/2}')}\right)$ by \eqref{eq:Holdestw}. \QED
%
%%%%%%%%%%%%%%%%%%%%%%%%%%%%%%%%%%%%%%%%%%%%%%%%%%%%%%%%%%%%%%%%%%%%%%%%%%%%%%%%%%%%%%%%%%%%%%%%%%%%

\noindent Let $V$  satisfy
\begin{equation}\label{eq:ipoV}
V\in C^1(\R^N\setminus\{0\}),
\quad |V(x)|+|x\cdot\nabla V(x)|\leq C|x|^{-2s} \text{ as } |x|\rightarrow 0
\text{ for some $C>0$.}
\end{equation}
Let $w\in H^1(B_R^+;t^{1-2s})$ solve
 \begin{equation}\label{eq:poho-wHextended}
\begin{cases}
  -  \dive(t^{1-2s}\nabla  w)+m^2t^{1-2s}w=0,&\text{in } B_{R}^+,\\
-\lim_{t\to 0^+}t^{1-2s}\frac{\partial w}{\partial
  t}(t,x)=\kappa_s V(x)w, &\text{on } B_R',
\end{cases}
\end{equation}
in a weak sense, i.e.,   for all  $\varphi\in C^ \infty_c(B_{R}^+\cup B_{R}')$,
\begin{equation}\label{eq:wH8}
{ \displaystyle\int_{\R^{N+1}_+}t^{1-2s}\nabla w\cdot\nabla \varphi\,dt\,dx}+m^2
  {\displaystyle\int_{\R^{N+1}_+}t^{1-2s}  w  \varphi\,dt\,dx}=
\kappa_s
\int_{B_R'}
V(x)w  \varphi\,dx.
\end{equation}
\noindent The following Pohozaev-type identity holds.
\begin{Theorem} \label{t:pohozaev}
Let $w$ be a solution to \eqref{eq:poho-wHextended} in sense of
\eqref{eq:wH8}, with $V$ satisfying \eqref{eq:ipoV}. Then,
for a.e. $r\in (0,R)$,
\begin{multline}\label{eq:poho}
  -\frac{N-2s}2 \int_{B_r^+} t^{1-2s}|\nabla w|^2dz
  -\frac{m^2(N+2-2s)}{2}\int_{B_r^+} t^{1-2s}  w^2dz\\+\frac{rm^2}{2}\int_{S_r^+} t^{1-2s}  w^2dS
  +\frac{r}{2}
\int_{S_r^+} t^{1-2s}|\nabla w|^2dS   \\
  =r\int_{S_r^+}t^{1-2s}\bigg|\frac{\partial w}{\partial
    \nu}\bigg|^2\,dS
  -\frac{\kappa_s}2\int_{B_r'}(NV(x)+\nabla V(x)\cdot x)w^2\,dx
 +
  \frac{r\kappa_s}2\int_{\partial B_r'}V(x)w^2\,dS'
\end{multline}
and
\begin{equation}\label{eq:poho2}
 \int_{B_r^+} t^{1-2s}|\nabla w|^2dz+m^2\int_{B_r^+} t^{1-2s} w^2dz
  =
\int_{S_r^+}t^{1-2s}\frac{\partial w}{\partial\nu}w\,dS+
\kappa_s
\int_{B'_r} V(x)w^2(x)  \,dx.
\end{equation}
\end{Theorem}

\proof
We have, on $B_R^+$, the formula
 \be\label{eq:diver} \dive\left(
  \frac{1}{2}t^{1-2s}|\n w|^2 z-t^{1-2s}(z\cdot \n w) \n w
\right)=\frac{N-2s}{2} t^{1-2s}|\n w|^2-(z\cdot \n w )
\dive(t^{1-2s}\n w).  \ee
Let $\rho<r<R$. Integrating by parts \eqref{eq:diver}
over the set
$$
O_\e:= (B_r^+\setminus{\ov{B_\rho^+}})\cap\{(t,x),\,t>\e\},
$$
with $\e>0$, we have
\begin{align*}
\frac{N-2s}{2}& \int_{O_\e} t^{1-2s}|\nabla w(z)|^2dz+m^2 \frac{N+2-2s}{2} \int_{O_\e} t^{1-2s}  w^2dz
-m^2\frac{r}{2}\int_{S_r^+\cap\{t>\e\}} t^{1-2s}  w^2dS\\
&+m^2\frac{\e^{2-2s}}{2}\int_{B_{\sqrt{r^2-\e^2}}'  \setminus
  B_{\sqrt{\rho^2-\e^2}}'}w^2(\e,x)dx
  +m^2\frac{\rho}{2}\int_{S_\rho^+\cap\{t>\e\}} t^{1-2s}  w^2dS
 \\
&=-\frac{1}{2}\e^{2-2s}\int_{B_{\sqrt{r^2-\e^2}}' \setminus
  B_{\sqrt{\rho^2-\e^2}}'}|\n w|^2(\e,x )dx\\
&+ \e^{2-2s}\int_{B_{\sqrt{r^2-\e^2}}'  \setminus
  B_{\sqrt{\rho^2-\e^2}}' }| w_t|^2(\e,x )dx \hspace{2cm}\\
&+\frac{r}{2}\int_{S_r^+\cap\{t>\e\}} t^{1-2s}|\nabla w|^2dS-
r\int_{S_r^+\cap\{t>\e\}}t^{1-2s}\bigg|\frac{\partial w}{\partial\nu}\bigg|^2dS\\
&-\frac{\rho}{2}\int_{S_\rho^+\cap\{t>\e\}} t^{1-2s}|\nabla w|^2dS
+\rho\int_{S_\rho^+\cap\{t>\e\}}t^{1-2s}\bigg|\frac{\partial w}{\partial\nu}\bigg|^2dS\\
&+ \int_{B_{\sqrt{r^2-\e^2}}'  \setminus
  B_{\sqrt{\rho^2-\e^2}}'}(x\cdot \n_x w(\e,x)) \,\e^{1-2s}w_t(\e,x)\,dx.
\end{align*}
We now claim that there exists a sequence $\e_n\to 0 $ such that
$$
\lim_{n\to \infty} \e_n^{2-2s}\left[   \int_{B_r' }|\n w|^2(\e_n,x )dx+
   \int_{B_r' }w^2(\e_n,x )dx \right]=0.
$$
If no such sequence exists, we  would have
$$
\liminf_{\e\to 0 }\e^{2-2s} \left[
\int_{B_r' }|\n w|^2(\e,x )dx +\int_{B_r' } w^2(\e,x )dx \right]\geq C>0
$$
and thus there exists $\e_0>0$ such that
$$
 \e^{2-2s}\bigg[\int_{B_r' }|\n w|^2(\e,x )dx
+\int_{B_r' } w^2(\e,x )dx\bigg]\geq \frac{C}{2}\quad \text{for all } \e\in(0,\e_0).
$$
It follows that, for all $\e\in(0,\e_0)$,
$$
 \frac{1}{2}\e^{1-2s}\int_{B_r' }|\n w|^2(\e,x )dx
+ \e^{1-2s}\int_{B_r' }w^2(\e,x )dx \geq \frac{C}{2\e}
$$
and so integrating the above inequality on $(0,\e_0)$ we contradict
the fact that $w\in H^1(B_{R}^+;t^{1-2s})$.

Next, from the Dominated Convergence Theorem, Lemma \ref{lem:reg}, and
Proposition \ref{p:Holdk},
we have  that
$$
\lim_{\e\to 0} \int_{B_{\sqrt{r^2-\e^2}}'  \setminus
  B_{\sqrt{\rho^2-\e^2}}'}(x\cdot \n_x w(\e,x))
\,\e^{1-2s}w_t(\e,x)\,dx= -\kappa_s\int_{B_r'\setminus B_\rho' }(x\cdot \n_x w)
\,V(x)w\,dx.
$$
We conclude (replacing $O_\e$ with $O_{\e_n}$, for a sequence $\e_n\to 0$) that
\begin{multline}\label{eq:pohozG}
\frac{N-2s}2 \int_{B_r^+\setminus{{B_\rho^+}}} t^{1-2s}|\nabla w(z)|^2dz
+m^2 \frac{N+2-2s}{2} \int_{B_r^+\setminus{{B_\rho^+}}} t^{1-2s}  w^2dz\\
-m^2\frac{r}{2}\int_{S_r^+ } t^{1-2s}  w^2dS +m^2\frac{\rho}{2}\int_{S_\rho^+ } t^{1-2s}  w^2dS
=
\frac{r}{2}\int_{S_r^+} t^{1-2s}|\nabla w|^2dS-r\int_{S_r^+}t^{1-2s}\bigg|\frac{\partial w}{\partial
    \nu}\bigg|^2\,dS\\
-\frac{\rho}{2}\int_{S_\rho^+} t^{1-2s}|\nabla w|^2dS-\rho\int_{S_\rho^+}t^{1-2s}\bigg|\frac{\partial w}{\partial
    \nu}\bigg|^2\,dS
- \kappa_s\int_{B_r'\setminus B_\rho' }(x\cdot \n_x w) \,V(x)w\,dx.
\end{multline}
Furthermore, integration by parts yields
\begin{align}\label{eq:9}
  \int_{B_r'\setminus B_\rho' }&(x\cdot \n_x w)
\,V(x)w\,dx=
-\frac{1}2\int_{B_r'\setminus B_\rho'}(NV(x)+\nabla V(x)\cdot x)w^2\,dx\\
\notag &+ \frac{r }2\int_{\partial B_r'}V(x)w^2\,{dS'}
-
  \frac{\rho }2\int_{\partial B_\rho'}V(x)w^2\,{dS'}.
\end{align}
 Since $w\in H^1(B_{R}^+;t^{1-2s})$, in view of Hardy  and Sobolev inequalities, there
 exists a sequence $\rho_n\to 0 $ such that
$$
\lim_{n\to \infty} \rho_n\bigg[\int_{S_{\rho_n}^+} t^{1-2s}[|\nabla
w|^2+w^2]dS + \int_{\partial
  B_{\rho_n}'}[|V(x)|+|x||\n V|]{w^2} dS'  \bigg] =0.
$$
Hence, taking $\rho=\rho_n$ and letting $n\to\infty$ in
\eqref{eq:pohozG} and \eqref{eq:9}, we obtain \eqref{eq:poho}.
Finally
\eqref{eq:poho2} follows the proof in \cite[Lemma 3.1]{FF}.
\QED

\section{The Almgren type frequency
  function}\label{sec:monot-prop}
In this section, we introduce  the \textit{Almgren
  frequency function} at the origin 0 for the extended problem
associated to the relativistic operator $\Dsm$ and
study its limit as $r\to 0^+$.
Let $R>0$
 and $w\in  H^1(B_{R}^+;t^{1-2s})$ be a nontrivial solution to
 \begin{equation}\label{eq:wHextended}
\begin{cases}
   - \dive(t^{1-2s}\nabla  w)+t^{1-2s}m^2w=0,&\text{in } B_{R}^+,\\
-\lim_{t\to 0^+}t^{1-2s}\frac{\partial w}{\partial
  t}(t,x)=\kappa_s
\Big(\frac{a(x/|x|)}{|x|^{2s}}w+hw\Big), &\text{on }B_R',
\end{cases}
\end{equation}
in the sense of \eqref{eq:wH8}.
Arguing as in \cite{FF}, it is easy to check that, for a.e. $r\in (0,R)$ and every $\widetilde\varphi\in
  C^\infty(\overline{B_r^+})$
\begin{equation}\label{eq:testbound}
\int_{B_r^+}t^{1-2s}\big(\nabla w\cdot\nabla\widetilde\varphi+m^2 w\widetilde\varphi\big)\,dz=
\int_{S_r^+}t^{1-2s}\frac{\partial w}{\partial\nu}\widetilde\varphi\,dS+
\kappa_s
\int_{B'_r}\bigg(
\frac{a(\frac{x}{|x|})}{|x|^{2s}}w+hw\bigg)\widetilde\varphi\,dx.
\end{equation}
The main result of this section is  the existence of  the limit as
$r\to 0^+$ of the \emph{Almgren's frequency} function (see
\cite{almgren} and \cite{GL}) associated to $w$
\begin{equation}\label{eq:1du}
  {\mathcal N}(r)=\frac{\displaystyle r \left[\int_{B_r^+}
t^{1-2s}\big(|\nabla w|^2+m^2w^2\big)\,dt\,dx-
\kappa_s
\int_{B_r'}\big(
\tfrac{a(x/|x|)}{|x|^{2s}}w^2+hw^2\big)\,dx\right]}{\displaystyle\int_{S_r^+}t^{1-2s}w^2 \, dS}.
\end{equation}
\begin{Theorem} \label{t:asym}
Let $w$ satisfy  \eqref{eq:wHextended}, with $s\in(0,1)$, $a\in
C^1({\mathbb S}^{N-1})$ satisfy
 \eqref{firsteig_strict_in},  and $h$ as in assumption
\eqref{eq:ipo1}. Then, letting
  ${\mathcal N}(r)$ as in (\ref{eq:1du}), there
there exists $k_0\in \N$, $k_0\geq1$, such that
\begin{equation}\label{eq:35du}
\lim_{r\to 0^+}{\mathcal
      N}(r)=
-\frac{N-2s}{2}+\sqrt{\bigg(\frac{N-2s}
    {2}\bigg)^{\!\!2}+\mu_{k_0}(a)}.
\end{equation}
Furthermore, if $\gamma$ denotes the limit in (\ref{eq:35du}), $M\geq
  1$ is the multiplicity of the eigenvalue
  $\mu_{j_0}(a)=\mu_{j_0+1}(a)=\cdots=\mu_{j_0+M-1}(a)$
  and
$\{\psi_i\}_{i=j_0}^{j_0+M-1}$ ($j_0\leq k_0\leq j_0+M-1$) is
  an $L^2({\mathbb
  S}^{N}_+;\theta_1^{1-2s})$-orthonormal basis for the eigenspace of
problem \eqref{eq:4} associated to $\mu_{k_0}(a)$, then
\begin{align*}
&\tau^{-\gamma}w(\tau\theta)\to
\sum_{i=j_0}^{j_0+M-1} \beta_i\psi_{i}(\theta)\quad \text{in }
C^{0,\alpha}(\overline{{\mathbb S}^{N}_+})  \quad \text{as }\tau\to 0^+,\\
& \tau^{-\gamma}w(0,\tau\theta')\to
\sum_{i=j_0}^{j_0+M-1} \beta_i\psi_{i}(0,\theta')\quad \text{in }
C^{1,\alpha}({\mathbb S}^{N-1}) \quad \text{as }\tau\to 0^+,
\end{align*}
for some $\alpha\in(0,1)$, where
\begin{align}
\label{eq:betai}  \beta_i&= R^{-\gamma}
\int_{{\mathbb S}^{N}_+}\theta_1^{1-2s}w(R\,\theta)
  \psi_i(\theta)\,dS
\\
\notag  &\quad-R^{-2\gamma-N+2s}\int_{0}^R\frac{\rho^{\gamma+N-1}}{2\gamma+N-2s}\bigg(
\kappa_s \int_{{\mathbb
      S}^{N-1}}h(\rho\theta')w(0,\rho\theta')\psi_i(0,\theta')\,dS'\\
\notag&\hskip6cm-m^2\rho^{2-2s}\int_{{\mathbb
  S}^{N}_+}\theta_1^{1-2s}w(\rho\theta)\psi_i(\theta)\,dS
\bigg) d\rho\\
 \notag &\quad +\int_{0}^R\frac{ \rho^{2s-\gamma-1}}{2\gamma+N-2s}\bigg(
\kappa_s \int_{{\mathbb
      S}^{N-1}}h(\rho\theta')w(0,\rho\theta') \psi_i(0,\theta')\,dS'\\
\notag&\hskip6cm
-m^2\rho^{2-2s}\int_{{\mathbb
  S}^{N}_+}\theta_1^{1-2s}w(\rho\theta)\psi_i(\theta)\,dS
 \bigg) d\rho
\end{align}
 for all $R>0$ such that $\overline{B_{R}'}\subset\Omega$
and $(\beta_{j_0},\beta_{j_0+1},\dots,\beta_{j_0+M-1})\neq(0,0,\dots,0)$.
\end{Theorem}

\noindent
From the Pohozaev-type identity \eqref{eq:poho} and
\eqref{eq:poho2} it follows
that,
for a.e. $r\in (0,R)$,
\begin{multline}\label{eq:2po}
  -\frac{N-2s}2\bigg[\int_{B_r^+} t^{1-2s}\big(|\nabla w|^2+m^2w^2\big)dz
  -\kappa_s
  \int_{B_r'}\frac{a(\frac{x}{|x|})}{|x|^{2s}}w^2 dx\bigg]\\
+\frac{r}{2}\bigg[
\int_{S_r^+} t^{1-2s}\big(|\nabla w|^2+m^2w^2\big)dS -\kappa_s
  \int_{\partial B_r'}\frac{a(\frac{x}{|x|})}{|x|^{2s}}w^2 dS'\bigg]\\
  =r\int_{S_r^+}t^{1-2s}\bigg|\frac{\partial w}{\partial
    \nu}\bigg|^2\,dS
  -\frac{\kappa_s}2\int_{B_r'}(Nh+\nabla h\cdot x)w^2\,dx
 +
  \frac{r\kappa_s}2\int_{\partial B_r'}hw^2\,dS'
+m^2\int_{B_r^+} t^{1-2s}w^2dz
\end{multline}
and
\begin{equation}\label{eq:3poho2}
  \int_{B_r^+}t^{1-2s}\big(|\nabla w|^2+m^2 w^2\big)\,dz
  -\kappa_s
  \int_{B'_r}\frac{a(\frac{x}{|x|})}{|x|^{2s}}w^2\,dz=
  \int_{S_r^+}t^{1-2s}\frac{\partial w}{\partial\nu}w\,dS+
  \kappa_s
  \int_{B'_r}
  hw^2\,dx.
\end{equation}
For every $r\in (0,R]$ we define
\begin{equation}\label{D(r)}
  D(r)=\frac{1}{r^{N-2s}} \bigg[
  \int_{B_r^+}
  t^{1-2s}\big(|\nabla w|^2+m^2w^2\big)\,dt\,dx-
  \kappa_s
  \int_{B_r'}\left(
    \frac{a(x/|x|)}{|x|^{2s}}w^2+h(x)w^2\right)\,dx\bigg]
\end{equation}
and
\begin{equation} \label{H(r)}
H(r)=\frac{1}{r^{N+1-2s}}\int_{S_r^+}t^{1-2s}w^2 \, dS
=
\int_{{\mathbb
      S}^{N}_+}\theta_1^{1-2s}w^2(r\theta)\,dS.
\end{equation}
\begin{Lemma} \label{l:hprime}
$H\in C^1(0,R)$ and
\begin{align}\label{H'}
 & H'(r)=\frac{2}{r^{N+1-2s}} \int_{S_r^+}t^{1-2s}w\frac{\partial
    w}{\partial\nu}
  \, dS, \quad \text{for every } r\in (0, R),\\
 &\label{H'2} H'(r)=\frac2r D(r), \quad \text{for every } r\in (0, R).
\end{align}
\end{Lemma}
\proof
The proof of \eqref{H'} can be performed arguing as in \cite[Lemma 3.8]{FF} and
using the regularity results of Lemma \ref{lem:reg} and Proposition
\ref{p:Holdk}.
The continuity of $H'$ on the interval $(0,R)$ follows by the
representation of $H'$ given above, Lemma \ref{lem:reg},
Propositions \ref{p:Hold}, \ref{p:Holdk}
and the
Dominated Convergence Theorem.
Finally, \eqref{H'2} follows from \eqref{H'}, \eqref{D(r)}, and \eqref{eq:3poho2}.
\QED

\noindent The regularity of the
function $D$ is established in the following lemma.

\begin{Lemma}\label{l:dprime}
  The function $D$ defined in (\ref{D(r)}) belongs to $W^{1,1}_{{\rm\
      loc}}(0, R)$ and
\begin{equation*}
  D'(r)=\frac{2}{r^{N+1-2s}} \bigg[ r\int_{S_r^+}
  t^{1-2s}\left|\frac{\partial w}{\partial \nu}\right|^2 dS
  -\kappa_s\int_{B_r'}\Big(sh+\frac12(\nabla h\cdot x)\Big)w^2\,dx
+m^2\int_{B_r^+} t^{1-2s}w^2dz
\bigg]
\end{equation*}
in a distributional sense and for a.e. $r\in (0,R)$.
\end{Lemma}

\proof
For any $r\in (0,r_0)$ let
\begin{equation}\label{I(r)}
I(r)=
\int_{B_r^+}
t^{1-2s}\big(|\nabla w|^2+m^2w^2\big)\,dt\,dx-
\kappa_s
\int_{B_r'}\left(
\frac{a(x/|x|)}{|x|^{2s}}w^2+h(x)w^2\right)\,dx
.
\end{equation}
Since $w\in  H^1(B_R^+;t^{1-2s})$, from \cite[Lemma 2.5]{FF}  we deduce that $I\in
 W^{1,1}(0,R)$ and
\begin{equation} \label{I'(r)}
I'(r) =
\int_{S_r^+}
t^{1-2s}\big(|\nabla w|^2+m^2w^2\big)\,dS-
\kappa_s
\int_{\partial B_r'}\left(
\frac{a(x/|x|)}{|x|^{2s}}w^2+h(x)w^2\right)\,dS'
\end{equation}
 for a.e. $r\in (0,R)$ and in the distributional sense.
Therefore $D\in W^{1,1}_{{\rm loc}}(0,R)$ and, using
\eqref{eq:2po}, (\ref{I(r)}), and (\ref{I'(r)}) into
\begin{align*}
  D'(r)=r^{2s-1-N}[-(N-2s)I(r)+rI'(r)],
\end{align*}
we obtain the conclusion.
\QED

\noindent
We prove now that, since $w\not\equiv 0$,
 $H(r)$ does not vanish for $r$ sufficiently small.

 \begin{Lemma} \label{welld}
There exists $R_0\in(0,R)$   such that
$H(r)>0$ for any $r\in (0,R_0)$, where $H$ is defined by (\ref{H(r)}).
\end{Lemma}

\proof
Let $R_0\in (0,R)$ such that $1-\kappa_s C_h R_0^\chi (C'_{a,N,s})^{-1}>0$, with $C'_{a,N,s}$ as in Corollary
\ref{c:hardyboundary_con2}. Suppose by contradiction that there exists $r_0\in(0,R_0)$ such that
$H(r_0)=0$. Then  $w= 0$ a.e. on $S_{r_0}^+$. From  \eqref{eq:3poho2} it follows that
\begin{equation*}
\int_{B_{r_0}^+}t^{1-2s}\big(|\nabla w|^2+m^2 w^2\big)\,dz
-\kappa_s
\int_{B'_{r_0}}\frac{a(\frac{x}{|x|})}{|x|^{2s}}w^2\,dz-
\kappa_s
\int_{B'_{r_0}}
hw^2\,dx=0.
\end{equation*}
From \eqref{eq:ipo1}, Corollaries \ref{c:hardyboundary_con} and
\ref{c:hardyboundary_con2}, it follows that
\begin{align*}
0&=
\int_{B_{r_0}^+}t^{1-2s}\big(|\nabla w|^2+m^2 w^2\big)\,dz
-\kappa_s
\int_{B'_{r_0}}\frac{a(\frac{x}{|x|})}{|x|^{2s}}w^2\,dz-
\kappa_s
\int_{B'_{r_0}}
hw^2\,dx\\
\notag & \geq
 C_{a,N,s} \Big(
1-\kappa_s C_h r_0^\chi (C'_{a,N,s})^{-1}\Big) \int_{B_{r_0}^+} t^{1-2s}|\nabla w|^2dz
,
\end{align*}
which, being $1-\kappa_s C_h r_0^\chi (C'_{a,N,s})^{-1} >0$, implies $w\equiv 0$ in
$B_{r_0}^+$ by Lemma \ref{l:hardyboundary}.  Classical
unique continuation principles for second order elliptic equations
with locally bounded coefficients (see e.g.  \cite{wolff}) allow to
conclude that $w=0$ a.e. in $B_R^+$, a contradiction.~\QED

\noindent Letting $R_0$ be as in Lemma \ref{welld} and  recalling \eqref{eq:1du},
the \emph{Almgren type frequency
  function}
\begin{equation}\label{N(r)}
{\mathcal N}(r)=\frac{D(r)}{H(r)}
\end{equation}
is well defined in $(0,R_0)$.

 \begin{Lemma}\label{mono} The function
   ${\mathcal N}$ defined in (\ref{N(r)}) belongs to $W^{1,1}_{{\rm
       loc}}(0, R_0)$ and
\begin{align}\label{formulona}
{\mathcal N}'(r)=\nu_1(r)+\nu_2(r)
\end{align}
in a distributional sense and for a.e. $r\in (0,R_0)$,
where
\begin{align}\label{eq:nu1}
\nu_1(r)=&\frac{2r\Big[
    \left(\int_{S_r^+}
  t^{1-2s}\left|\frac{\partial w}{\partial \nu}\right|^2 dS\right)
    \left(\int_{S_r^+}
  t^{1-2s}w^2\,dS\right)-\left(
\int_{S_r^+}
  t^{1-2s}w\frac{\partial w}{\partial \nu}\, dS\right)^{\!2} \Big]}
{\left(
\int_{S_r^+}
  t^{1-2s}w^2\,dS\right)^2}
\end{align}
and
\begin{equation}\label{eq:nu2}
  \nu_2(r)= \frac{ 2m^2\int_{B_r^+} t^{1-2s}w^2dz-\kappa_s\int_{B_r'} (2sh+\nabla h\cdot
    x)w^2\,dx}{\int_{S_r^+} t^{1-2s}w^2\,dS}.
\end{equation}
\end{Lemma}
\proof It follows from Lemmas \ref{l:hprime}, \ref{welld}, and
  \ref{l:dprime}.
\QED

\begin{Lemma} \label{l:stimasotto}
 Let   ${\mathcal N}$ be the function defined in (\ref{N(r)}).
There
exist $\tilde R\in (0,R_0)$ and a  constant
$\overline{C}>0$  such that
\begin{multline}\label{eq:13}
\int_{B_r^+}
t^{1-2s}\big(|\nabla w|^2+m^2w^2\big)\,dt\,dx-
\kappa_s
\int_{B_r'}\left(
\frac{a(x/|x|)}{|x|^{2s}}w^2+h(x)w^2\right)\,dx\\
\geq
-\bigg(\frac{N-2s}{2r}\bigg)\int_{S_r^+}t^{1-2s}w^2dS +\overline{C}\bigg(
\int_{B_r'}\dfrac{w^2}{|x|^{2s}}\,dx+\int_{B_r^+}
t^{1-2s}|\nabla w|^2\,dt\,dx+
\int_{B_r^+}
t^{1-2s}\frac{w^2}{|z|^2}\,dt\,dx\bigg)
\end{multline}
and
\begin{equation}\label{Nbelow}
   {\mathcal N}(r)>-\frac{N-2s}{2}
 \end{equation}
for every $r\in(0,\tilde R)$.
\end{Lemma}
\proof
From Lemma \ref{l:hardyboundary} and  Corollary
\ref{c:hardyboundary_con2}, it follows that
\begin{multline*}
\int_{B_r^+}
t^{1-2s}\big(|\nabla w|^2+m^2w^2\big)\,dt\,dx-
\kappa_s
\int_{B_r'}\left(
\frac{a(x/|x|)}{|x|^{2s}}w^2+h(x)w^2\right)\,dx
+\bigg(\frac{N-2s}{2r}\bigg)\int_{S_r^+}t^{1-2s}w^2dS\\
\geq \Big(
1-\frac{m^2r^2}{\mu_1(a)+(\frac{N-2s}{2})^2}-\frac{\kappa_sC_h r^\chi}{C'_{N,a,s}}
\Big)
\bigg(\int_{B_r^+}
t^{1-2s}|\nabla w|^2\,dt\,dx\\
-\kappa_s
\int_{B_r'}
\frac{a(x/|x|)}{|x|^{2s}}w^2\,dx+\bigg(\frac{N-2s}{2r}\bigg)\int_{S_r^+}t^{1-2s}w^2dS\bigg)
\end{multline*}
for every $r\in(0,R_0)$. The conclusion follows from the above
estimate, choosing $r$ sufficiently small and using Lemma \ref{l:hardyboundary} and  Corollaries
\ref{c:hardyboundary_con},
\ref{c:hardyboundary_con2}.
\QED

\begin{Lemma} \label{l:stimanu2}
 Let $\tilde R$ be as in Lemma \ref{l:stimasotto} and
  $\nu_2$ as in (\ref{eq:nu2}).
There
exists
  $C_1>0$ such that
\begin{equation*}
|\nu_2(r)|\leq C_1\left[{\mathcal
      N}(r)+\frac{N-2s}{2}\right]r^{-1+\chi}
\end{equation*}
for a.e. $r\in (0,\tilde R)$.
\end{Lemma}

\proof
From \eqref{eq:ipo1} and  (\ref{eq:13}) we deduce that
\begin{align*}
\left|
\int_{B_r'} (2sh(x)+\nabla h(x)\cdot
    x)w^2\,dx\right|&\leq
2C_h r^\chi
\int_{B_r'} \frac{|w|^2}{|x|^{2s}}\,dx
\leq  2C_h\overline{C}^{-1}\,
r^{\chi+N-2s}\left[D(r)+{\textstyle{\frac{N-2s}{2}}}H(r)\right],
\end{align*}
and, therefore, for any $r\in (0,\tilde R)$, we have that
\begin{align}\label{B00}
  \left|
\frac{\int_{B_r'} (2sh(x)+\nabla h(x)\cdot
    x)w^2\,dx}{\int_{S_r^+} t^{1-2s}w^2\,dS}\right|&\leq 2C_h\overline{C}^{-1}
  \,r^{-1+\chi} \frac{D(r)+\frac{N-2s}{2}H(r)}{H(r)} \\
&\notag=
  2C_h\overline{C}^{-1}\, r^{-1+\chi}\left[{\mathcal N}(r)+
    \frac{N-2s}{2}\right].
\end{align}
On the other hand, from (\ref{eq:13}) it also follows that
\begin{equation}\label{eq:stm}
 \left|
\frac{ \int_{B_r^+} t^{1-2s}w^2dz}
{\int_{S_r^+} t^{1-2s}w^2\,dS}\right|
\leq \overline{C}^{-1} r \left[{\mathcal N}(r)+
    \frac{N-2s}{2}\right].
\end{equation}
Combining \eqref{B00} with \eqref{eq:stm} we obtain the stated
estimate.
\QED

\begin{Lemma} \label{l:stima_N_sopra} Let $\tilde R$ be as in Lemma \ref{l:stimasotto},
   ${\mathcal N}$ as in (\ref{N(r)}) and $H$ as in \eqref{H(r)}.Then
\begin{enumerate}[\rm (i)]
\item there exist a positive
  constant $C_2>0$   such that
${\mathcal N}(r)\leq C_2$
for all $r\in (0,\tilde R)$;
\item the limit
$\gamma:=\lim_{r\rightarrow 0^+} {\mathcal N}(r)$
exists and is finite;
\item   there exists a constant $K_1>0$ such that
$H(r)\leq K_1 r^{2\gamma}$ for all $r\in (0,\tilde R)$;
\item for any $\sigma>0$ there exists a constant
$K_2(\sigma)>0$ depending on $\sigma$ such that
$H(r)\geq K_2(\sigma)\,
  r^{2\gamma+\sigma}$ for all $r\in (0,\tilde R)$.
\end{enumerate}
\end{Lemma}
\proof
By Lemma \ref{mono}, Schwarz's inequality, and Lemma \ref{l:stimanu2}, we
 obtain
\begin{equation}\label{eq:40}
\bigg({\mathcal N}+\frac {N-2s}2\bigg)'(r)\geq
\nu_2(r)\geq
 -C_1\left[{\mathcal
      N}(r)+\frac{N-2s}{2}\right]r^{-1+\chi}
\end{equation}
for a.e. $r\in (0,\tilde R)$.
Integration over $(r,\tilde R)$ yields
\begin{equation*}
{\mathcal N}(r)\leq -\frac {N-2s}2+
\left({\mathcal
      N}(\tilde R)+\frac{N-2s}{2}\right)
e^{\frac{C_1}{\chi}\tilde R^\chi}
\end{equation*}
for any $r\in (0,\tilde R)$, thus proving claim (i).

  By Lemmas \ref{l:stimanu2} and \ref{l:stima_N_sopra}, the function
  $\nu_2$ defined in (\ref{eq:nu2}) belongs to $L^1(0,\tilde R)$.
  Hence, by Lemma \ref{mono} and Schwarz's inequality, ${\mathcal N}'$
  is the sum of a nonnegative function and of a $L^1$-function on
  $(0,\tilde R)$.  Therefore
$$
{\mathcal N}(r)={\mathcal N}(\tilde R)-\int_r^{\tilde R} {\mathcal N}'(\rho)\, d\rho
$$
admits a limit as $r\rightarrow 0^+$ which is necessarily finite in view of
(\ref{Nbelow}) and part (i). Claim (ii) is thereby proved.

  By (ii) ${\mathcal N}'\in L^1(0,\tilde R)$ and, by (i), ${\mathcal N}$ is bounded, then from
  (\ref{eq:40}) and (i) it follows that
  \begin{equation*}
{\mathcal N}(r)-\gamma=\int_0^r
    {\mathcal N}'(\rho) \, d\rho\geq -C_3 r^\chi
\end{equation*}
for all $r\in(0, \tilde R)$.
Therefore by (\ref{H'2}) and \eqref{N(r)},
we deduce that, for all $r\in(0,\tilde R)$,
$$
\frac{H'(r)}{H(r)}=\frac{2\,{\mathcal N}(r)}{r}\geq
\frac{2\gamma}{r}-2C_3 r^{-1+\chi},
$$
which, after integration over the interval $(r, \tilde R)$, yields (iii).

From (ii) it follows that, for any $\sigma>0$ there exists $r_\sigma>0$ such
that ${\mathcal N}(r)<\gamma+\frac\sigma2$ for any $r\in (0,r_\sigma)$ and
hence
$$
\frac{H'(r)}{H(r)}=\frac{2\,{\mathcal N}(r)}{r}<\frac{2\gamma+\sigma}{r}
\quad \text{for all } r\in (0,r_\sigma).
$$
Integrating over the interval $(r,r_\sigma)$ and by continuity of $H$
outside $0$, we obtain~(iv).~\QED

\subsection{The blow-up argument  }\label{sec:blow-up}

\begin{Lemma}\label{l:blowup}
Let $w$ satisfy  \eqref{eq:wHextended}, with $s\in(0,1)$, $h$ as in
assumption \eqref{eq:ipo1} and
$a\in C^1({\mathbb S}^{N-1})$ satisfy \eqref{firsteig_strict_in}. Let $\gamma:=\lim_{r\rightarrow 0^+} {\mathcal
    N}(r)$ as in Lemma \ref{l:stima_N_sopra}. Then
\begin{itemize}
\item[\rm (i)] there exists $k_0\in \N$, $k_0\geq1$, such that
  $\gamma=-\frac{N-2s}{2}+\sqrt{\big(\frac{N-2s}
    {2}\big)^{2}+\mu_{k_0}(a)}$;
\item[\rm (ii)] for every sequence $\tau_n\to0^+$, there exist a subsequence
$\{\tau_{n_k}\}_{k\in\N}$ and an eigenfunction $\psi$ of problem
\eqref{eq:4} associated to the eigenvalue $\mu_{k_0}(a)$ such that
$\|\psi\|_{L^2({\mathbb
  S}^{N}_+;\theta_1^{1-2s})}=1$ and
\[
\frac{w(\tau_{n_k}z)}{\sqrt{H(\tau_{n_k})}}\to
|z|^{\gamma}\psi\Big(\frac z{|z|}\Big)
\]
 strongly in $H^1(B_r^+;t^{1-2s})$ and in  $C^{0,\alpha}_{\rm
    loc}(\overline{B_r^+}\setminus\{0\})$  for some $\alpha\in (0,1)$
  and all $r\in(0,1)$
and
\[
\frac{w(0,\tau_{n_k}x)}{\sqrt{H(\tau_{n_k})}}\to
|x|^{\gamma}\psi\Big(0,\frac x{|x|}\Big)
\]
in  $C^{1,\alpha}_{\rm
    loc}(B_1'\setminus\{0\})$.
\end{itemize}
\end{Lemma}
\proof
Let us set
\begin{equation}\label{eq:wtau}
w^\tau(z)=\frac{w(\tau z)}{\sqrt{H(\tau)}}.
\end{equation}
We notice that $\int_{S_1^+}t^{1-2s}|w^{\tau}|^2dS=1$. Moreover, by
scaling and Lemma \ref{l:stima_N_sopra}, part~(i),
\begin{equation*}
  \int_{B_{1}^+} t^{1-2s}\Big(|\nabla w^\tau(z)|^2
+m^2\tau^2 | w^\tau(z)|^2\Big) dz
 -
\kappa_s
\int_{B_1'}\bigg(
\dfrac{a(\frac x{|x|})}{|x|^{2s}}+\tau^{2s}h(\tau
x)\bigg)|w^\tau |^2dx
  ={\mathcal N}(\tau)\leq C_2
\end{equation*}
for every $\tau\in(0,\tilde R)$, whereas, from \eqref{eq:13},
\begin{multline*}
  \mathcal N(\tau)\geq \frac{\tau^{-N+2s}}{H(\tau)}\bigg(
-\bigg(\frac{N-2s}{2\tau}\bigg)\int_{S_\tau^+}t^{1-2s}w^2dS
+\overline{C}\int_{B_\tau^+}
t^{1-2s}|\nabla w|^2\,dt\,dx\bigg)\\
=-\frac{N-2s}{2}+\overline{C} \int_{B_{1}^+} t^{1-2s}|\nabla w^\tau(z)|^2 dz
\end{multline*}
for every $\tau\in(0,\tilde R)$. From the above estimates,
 $\{w^\tau\}_{\tau\in(0,\tilde R)}$ is bounded in $H^1(B_1^+;t^{1-2s})$.
Therefore, for any given sequence $\tau_n\to 0^+$, there exists a
subsequence $\tau_{n_k}\to0^+$ such that $w^{\tau _{n_k}}\weakly
\widetilde w$ weakly in $H^1(B_1^+;t^{1-2s})$ for some $\widetilde
w\in H^1(B_1^+;t^{1-2s})$.
Moreover,  $\int_{S_1^+}t^{1-2s}|\widetilde w|^2dS=1$ due to compactness of the trace embedding $H^1(B_1^+;t^{1-2s})\hookrightarrow \hookrightarrow L^2(S_1^+;t_1^{1-2s})$. In particular $\widetilde w\not\equiv
0$.

For every small  $\tau\in (0,\tilde R)$, $w^\tau$ satisfies
\begin{equation} \label{eqlam}
\begin{cases}
   - \dive(t^{1-2s}\nabla  w^\tau)+\tau^2 t^{1-2s} m^2w^{\tau}=0,&\text{in }B_1^+,\\
-\lim_{t\to 0^+}t^{1-2s}\frac{\partial w^\tau}{\partial
  t}=\kappa_s
\Big(\frac{a(x/|x|)}{|x|^{2s}}w^\tau+\tau^{2s}h(\tau
x)w^\tau\Big), &\text{on }B_1',
\end{cases}
\end{equation}
in a weak sense, i.e.
\begin{equation*}
\int_{B_1^+}t^{1-2s}\big(\nabla w^\tau\!\cdot\!\nabla\widetilde\varphi+m^2\tau^2 w^\tau\widetilde\varphi\big)\,dz=
\kappa_s
\int_{B'_1}\bigg(
\frac{a(\frac{x}{|x|})}{|x|^{2s}}w^\tau+\tau^{2s}h(\tau
x)w^\tau \bigg)\widetilde\varphi(0,x)\,dx
\end{equation*}
for all $\widetilde\varphi \in
H^1(B_1^+;t^{1-2s})$ s.t. $ \widetilde\varphi=0$ on $\S^N_+$ and,  for such  $\widetilde\varphi$,
 by \eqref{eq:ipo1} and \cite[Lemma 2.5]{FF},
\[
 \tau^{2s}\int_{B_1'}h(\tau
x)w^\tau\widetilde\varphi(0,x)\,dx=o(1)\quad\text{as }\tau\to 0^+
\]
and, by \cite[Lemma 2.4]{FF},
\[
\tau^2\int_{B_1^+}t^{1-2s} w^\tau\widetilde\varphi\,dz=o(1)\quad\text{as }\tau\to 0^+.
\]
From weak convergence
 $w^{\tau _{n_k}}\weakly
\widetilde w$ in $H^1(B_1^+;t^{1-2s})$, we can pass to
the limit in \eqref{eqlam} along the sequence $\tau _{n_k}$
and obtain  that $\widetilde w$ weakly solves
\begin{equation}\label{eq:extended_limit}
\begin{cases}
    \dive(t^{1-2s}\nabla \widetilde w)=0,&\text{in }B_1^+,\\
-\lim_{t\to 0^+}t^{1-2s}\frac{\partial \widetilde w}{\partial
  t}=\kappa_s
\frac{a(x/|x|)}{|x|^{2s}}\widetilde w, &\text{on }B_1'.
\end{cases}
\end{equation}
From Proposition \ref{p:Hold}, we have that
\begin{equation*}
  w^{\tau _{n_k}}\to
  \widetilde w
\quad\text{in }C^{0,\alpha}_{\rm loc}(\overline{B_r^+}\setminus\{0\}),
\end{equation*}
while Proposition \ref{p:Holdk} and Lemma \ref{lem:reg} imply that
\begin{equation}\label{eq:17}
  \nabla_xw^{\tau _{n_k}}\to
  \nabla_x\widetilde w,\quad\text{and}\quad
  t^{1-2s}\frac{\partial w^{\tau _{n_k}}}{\partial t}\to
  t^{1-2s}\frac{\partial \widetilde w}{\partial t}
  \quad\text{in }C^{0,\alpha}_{\rm loc}(\overline{B_r^+}\setminus\{0\})
\end{equation}
for some $\alpha\in (0,1)$ and all $r\in(0,1)$.
 By \eqref{eq:ipo1}, \cite[Lemma 2.5]{FF}, and boundedness of
$\{w^\tau\}_{\tau\in(0,\tilde R)}$ in $H^1(B_1^+;t^{1-2s})$, it
follows that
\begin{equation}\label{eq:conze1}
 \tau^{2s}\int_{B_1'}h(\tau
x)|w^\tau|^2\,dx=o(1)\quad\text{as }\tau\to 0^+,
\end{equation}
and, by \cite[Lemma 2.4]{FF},
\begin{equation}\label{eq:conze2}
\tau^2\int_{B_1^+}t^{1-2s} |w^\tau|^2\,dz=o(1)\quad\text{as }\tau\to 0^+.
\end{equation}
Multiplying  equation \eqref{eqlam} with $w^\tau$, integrating in
$B_r^+$, and using
\eqref{eq:17}, \eqref{eq:conze1}, \eqref{eq:conze2}, and Corollary \ref{c:hardyboundary_con}, we easily obtain that $\|w^{\tau
  _{n_k}}\|_{H^1(B_r^+;t^{1-2s})}\to \|\widetilde
w\|_{H^1(B_r^+;t^{1-2s})}$ for all $r\in(0,1)$, and hence
\begin{equation}\label{strongH1}
w^{\tau _{n_k}}\to \widetilde w\quad\text{in }H^1(B_r^+;t^{1-2s})
\end{equation}
for any $r\in (0,1)$.
For any $r\in (0,1)$ and $k\in \N$, let us define
the functions
\begin{multline*}
D_k(r)=\frac{1}{r^{N-2s}} \bigg[\int_{B_r^+}
t^{1-2s}\big(|\nabla w^{\tau_{n_k}}|^2+m^2 \tau_{n_k}^2|w^{\tau_{n_k}}|^2\big)
\,dt\,dx\\
 -
\kappa_s
\int_{B_r'}\bigg(
\dfrac{a(\frac x{|x|})}{|x|^{2s}}+\tau_{n_k}^{2s}h(\tau_{n_k}
x)\bigg)|w^{\tau_{n_k}} |^2dx\bigg]
\end{multline*}
and
\begin{equation*}
H_k(r)=\frac{1}{r^{N+1-2s}}\int_{S_r^+}t^{1-2s}|w^{\tau_{n_k}}|^2 \, dS.
\end{equation*}
Direct calculations show that ${\mathcal
    N}_k(r):=\frac{D_k(r)}{H_k(r)}=\frac{D(\tau_{n_k}r)}{H(\tau_{n_k}r)}
  ={\mathcal N}(\tau_{n_k}r)$ for all $r\in (0,1)$.
From (\ref{strongH1}), \eqref{eq:conze1}, and \eqref{eq:conze2}, it follows that,
for any fixed $r\in (0,1)$, $D_k(r)\to \widetilde D(r)$,
where
\begin{equation*}
  \widetilde D (r)=
\frac{1}{r^{N-2s}} \bigg[\int_{B_r^+}
t^{1-2s}|\nabla \widetilde w|^2\,dt\,dx-
\kappa_s
\int_{B_r'}
\dfrac{a(\frac{x}{|x|})}{|x|^{2s}}\widetilde w^2\,dx\bigg]
\quad \text{for all } r\in (0,1).
\end{equation*}
The compactness of the trace embedding
 $H^1(B_r^+;t^{1-2s})\hookrightarrow \hookrightarrow
 L^2(S_r^+;t_1^{1-2s})$ ensures that,  for
 every  $r\in (0,1)$, $H_k(r)\to \widetilde H(r)$,
where
\begin{equation*}
\widetilde  H(r)=\frac{1}{r^{N+1-2s}}\int_{S_r^+}t^{1-2s}\widetilde
w^2 \, dS.
\end{equation*}
Arguing as in Lemma \ref{welld}, we can easily prove that $\widetilde H(r)>0$ for all
$r\in(0,1)$ and the function
\begin{equation}\label{eq:deftilden}
\widetilde {\mathcal N}(r):=\frac{\widetilde  D(r)}{\widetilde  H(r)}
\end{equation}
is well defined for $r\in (0,1)$. From
Lemma \ref{l:stima_N_sopra} part (ii), we deduce  that
\begin{equation*}
\widetilde {\mathcal
    N}(r)=\lim_{k\to \infty} {\mathcal
    N}(\t_{n_k}r)=\gamma
\end{equation*}
for all $r\in (0,1)$.
In particular, $\widetilde {\mathcal N}$ is  constant in $(0,1)$ and hence $\widetilde {\mathcal
  N}'(r)=0$ for any $r\in (0,1)$.  By \eqref{eq:extended_limit} and Lemma
\ref{mono} with $h\equiv 0$ and $m=0$, we obtain
\begin{equation*}
    \left(\int_{S_r^+}
  t^{1-2s}\left|\frac{\partial \widetilde w}{\partial \nu}\right|^2 dS\right) \cdot
    \left(\int_{S_r^+}
  t^{1-2s}\widetilde w^2\,dS\right)-\left(
\int_{S_r^+}
  t^{1-2s}\widetilde w\frac{\partial \widetilde w}{\partial \nu}\,
  dS\right)^{\!\!2}=0
\end{equation*}
for all $r\in (0,1)$,
which  implies  that $\widetilde w$ and $\frac{\partial \widetilde  w}{\partial \nu}$
have the same direction as vectors in $L^2(S_r^+;t^{1-2s})$ and hence
there exists a function $\eta=\eta(r)$ such that
$\frac{\partial \widetilde w}{\partial \nu}(r,\theta)=\eta(r)
\widetilde w(r,\theta)$ for all $r\in(0,1)$ and $\theta\in {\mathbb S}^N_+$.
After integration we obtain
\begin{equation} \label{separate}
\widetilde w(r,\theta)=e^{\int_1^r \eta(s)ds} \widetilde w(1,\theta)
=\varphi(r) \psi(\theta), \quad  r\in(0,1), \ \theta\in {\mathbb S}^N_+,
\end{equation}
where $\varphi(r)=e^{\int_1^r \eta(s)ds}$ and $\psi(\theta)=\widetilde w(1,\theta)$.
From \eqref{eq:extended_limit} and  (\ref{separate}), it follows that
$$
\begin{cases}
\frac1{r^N}\big(r^{N+1-2s}\varphi'\big)'\theta_1^{1-2s}\psi(\theta)+
r^{-1-2s}\varphi(r)\dive\nolimits_{{\mathbb
    S}^{N}}(\theta_1^{1-2s}\nabla_{{\mathbb S}^{N}}\psi(\theta))
=0,\\
-\lim_{\theta_1\to 0^+} \theta_1^{1-2s}\nabla_{{\mathbb
    S}^{N}}\psi(\theta)\cdot {\mathbf
  e}_1=\kappa_s a(\theta')\psi(0,\theta').
\end{cases}
$$
Taking $r$ fixed we deduce that $\psi$ is an eigenfunction of the
eigenvalue problem \eqref{eq:4}. If $\mu_{k_0}(a)$ is the corresponding
eigenvalue then $\varphi(r)$ solves the equation
$$
\varphi''(r)+\frac{N+1-2s}r\varphi'-\frac{\mu_{k_0}(a)}{r^2}\varphi(r)=0
$$
and hence $\varphi(r)$ is of the form
$$
\varphi(r)=c_1 r^{\sigma_{k_0}^+}+c_2 r^{\sigma_{k_0}^-}
$$
for some $c_1,c_2\in\R$, where
\begin{equation*}
  \sigma^+_{k_0}=-\tfrac{N-2s}{2}+\sqrt{\big(\tfrac{N-2s}
    {2}\big)^{2}+\mu_{k_0}(a)}\quad\text{and}\quad
  \sigma^-_{k_0}=-\tfrac{N-2s}{2}-\sqrt{\big(\tfrac{N-2s}{2}
    \big)^{2}+\mu_{k_0}(a)}.
\end{equation*}
Since the function
$|x|^{\sigma_{k_0}^-}\psi(\frac{x}{|x|})\notin
L^2(B_1';|x|^{-2s})$ and hence $|z|^{\sigma_{k_0}^-}\psi(\frac{z}{|z|})
\notin  H^1(B_1^+;t^{1-2s})$  in virtue of Lemma \cite[Lemma 2.5]{FF}, we deduce that $c_2=0$ and  $\varphi(r)=c_1
r^{\sigma_{k_0}^+}$. Moreover, from $\varphi(1)=1$, we obtain that $c_1=1$ and
then
\begin{equation} \label{expw}
\widetilde w(r,\theta)=r^{\sigma_{k_0}^+} \psi(\theta),  \quad
\text{for all }r\in (0,1)\text{ and }\theta\in {\mathbb S}^N_+.
\end{equation}
Substituting \eqref{expw} into \eqref{eq:deftilden}, we obtain that $\gamma=\widetilde{\mathcal N}(r)=\frac{\widetilde D(r)}{\widetilde H(r)}=\sigma_{k_0}^+$.
This completes the proof of the lemma.
\QED
\begin{Lemma} \label{l:limite}
If $w$ satisfies \eqref{eq:wHextended}, $H$ is defined in
\eqref{H(r)}, and  $\gamma:=\lim_{r\rightarrow 0^+} {\mathcal
    N}(r)$ is as in Lemma \ref{l:stima_N_sopra}, then the limit
$\lim_{r\to0^+}r^{-2\gamma}H(r)$
exists and it is finite.
\end{Lemma}
\proof
In view of  Lemma \ref{l:stima_N_sopra}, part (iii), it is sufficient to prove that the limit
exists. From (\ref{H(r)}), (\ref{H'2}), and Lemma \ref{l:stima_N_sopra} it
follows that
\begin{equation*}
\frac{d}{dr} \frac{H(r)}{r^{2\gamma}} =2r^{-2\gamma-1} (D(r)-\gamma
H(r))=2r^{-2\gamma-1} H(r) \int_0^r {\mathcal N}'(\rho) d\rho,
\end{equation*}
which, by integration over $(r,\tilde R)$, yields
\begin{equation}\label{inte}
  \frac{H(\tilde R)}{\tilde R^{2\gamma}}-
  \frac{H(r)}{r^{2\gamma}}=\int_r^{\tilde R}
f_1(\rho)\,d\rho
+\int_r^{\tilde R}
f_2(\rho) d\rho
\end{equation}
where $f_i(\rho)=2\rho^{-2\gamma-1}
  H(\rho) \left( \int_0^\rho \nu_i(t) dt \right)$, $i=1,2$, and
$\nu_1$ and $\nu_2$ are as in (\ref{eq:nu1}) and (\ref{eq:nu2}).
Since, by Schwarz's inequality, $\nu_1\geq 0$, we have that
$\lim_{r\to 0^+} \int_r^{\tilde R} f_1(\rho)d\rho$
exists.  On the other hand, by Lemmas \ref{l:stimanu2} and
\ref{l:stima_N_sopra}, we have that
\begin{equation*}
  \left| f_2(\rho)\right|\leq \frac{2K_1C_1}{\chi}\Big(C_2+\frac {N-2s}2\Big)\rho^{-1+\chi}
\end{equation*}
for all $\rho\in(0,\widetilde R)$, which proves that $f_2\in L^1(0,\widetilde R)$.
Hence both terms at the right hand side of
(\ref{inte}) admit a limit as $r\to 0^+$ thus completing the proof.
\QED

\noindent From Lemma \ref{l:blowup}, the following point-wise estimate for
solutions to \eqref{eq:wHextended} follow.
\begin{Lemma}\label{l:stima0}
 Let  $w$
satisfy \eqref{eq:wHextended}. Then there exists $C_4>0$ and $\bar
r\in(0,\tilde R)$ such that $|w(z)|\leq C_4|z|^\gamma$ for all $z\in
\overline{B_{\bar r}^+}$, where $\gamma:=\lim_{r\rightarrow 0^+}
{\mathcal
    N}(r)$ is as in Lemma \ref{l:stima_N_sopra}.
\end{Lemma}
\proof
We first claim that
\begin{equation}
  \label{eq:ilemm}
  \sup_{S_{r}^+}|w|^2=O(H(r))\quad\text{as }r\to 0^+.
\end{equation}
 In order to prove \eqref{eq:ilemm}, we argue by contradiction and
assume that there exists a sequence $\tau_n\to 0^+$ such that
\begin{equation*}
  \sup_{\theta\in{\mathbb S}^N_+}\Big|w\Big(\frac{\tau_n}2\theta\Big)\Big|^2>n H \Big (\frac{\tau_n}2 \Big),
\end{equation*}
i.e., defining $w^\tau$ as in \eqref{eq:wtau},
\begin{equation}\label{eq:18}
  \sup_{x\in S_{1/2}^+}|w^{\tau_n}(z)|^2>2^{N+1-2s}n\int_{S_{1/2}^+}t^{1-2s}|w^{\tau_n}(z)|^2dS.
\end{equation}
From Lemma \ref{l:blowup}, along a subsequence $\tau_{n_k}$ we have
that $w^{\tau_{n_k}}\to
|z|^{\gamma}\psi\big(\frac z{|z|}\big)$ in  $C^{0,\alpha}_{\rm
    loc}(\overline{S_{1/2}^+})$, for some $\psi$ eigenfunction of problem
\eqref{eq:4}, hence passing to the limit in \eqref{eq:18} gives rise
to a contradiction and claim \eqref{eq:ilemm} is proved.
The conclusion follows from combination of \eqref{eq:ilemm} and part
(iii) of Lemma \ref{l:stima_N_sopra}.~\QED

\noindent We will now prove that
$\lim_{r\to 0^+} r^{-2\gamma} H(r)$ is strictly positive.

\begin{Lemma} \label{l:limitepositivo} Under the same assumptions of
  Lemma \ref{l:limite}, $\lim_{r\to0^+}r^{-2\gamma}H(r)>0$.
\end{Lemma}
\proof
 For all $k\geq1$, let $\psi_k$ be as in \eqref{angular}, i.e. $\psi_k$
is a $L^2({\mathbb
  S}^{N}_+;\theta_1^{1-2s})$-normalized eigenfunction of problem
\eqref{eq:4} associated to the eigenvalue $\mu_k(a)$ and
 $\{\psi_k\}_k$ is an orthonormal basis of $L^2({\mathbb
  S}^{N}_+;\theta_1^{1-2s}) $.
From Lemma \ref{l:blowup}
there exist
$j_0,M\in\N\setminus\{0\}$, such that $M$ is the multiplicity of the
eigenvalue $\mu_{j_0}(a)=\mu_{j_0+1}(a)=\cdots=\mu_{j_0+M-1}(a)$ and
\begin{equation*}
  \gamma=\lim_{r\rightarrow 0^+} {\mathcal N}(r)=-\frac{N-2s}{2}+\sqrt{\bigg(\frac{N-2s}
    {2}\bigg)^{\!\!2}+\mu_{i}(a)},
  \quad i=j_0,\dots,j_0+M-1.
\end{equation*}
Let us expand $w$ as $w(z)=w(\tau
\theta)=\sum_{k=1}^\infty\varphi_k(\tau)\psi_k(\theta)$,
where $\tau=|z|\in(0,R]$, $\theta=z/|z|\in{{\mathbb S}^{N}_+}$, and
\begin{equation*}
  \varphi_k(\tau)=\int_{{\mathbb S}^{N}_+}\theta_1^{1-2s}w(\tau\,\theta)
  \psi_k(\theta)\,dS.
  \end{equation*}
The Parseval identity yields
\begin{equation}\label{eq:17bis}
H(\tau)=\int_{{\mathbb
    S}^{N}_+}\theta_1^{1-2s}w^2(\tau\theta)\,dS=
\sum_{k=1}^{\infty}\varphi_k^2(\tau),\quad\text{for all }0<\tau\leq R.
\end{equation}
In particular, from Lemma \ref{l:stima_N_sopra} (iii) and \eqref{eq:17bis} it follows that, for all $k\geq1$,
\begin{equation}\label{eq:23}
 \varphi_k(\tau)=O(\tau^{\gamma})\quad\text{as
 }\tau\to0^+.
\end{equation}
Equations \eqref{eq:wHextended} and \eqref{angular} imply that, for every $k$,
\begin{equation*}
-\varphi_k''(\tau)-\frac{N+1-2s}{\tau}\varphi_k^\prime(\tau)+
\frac{\mu_k(a)}{\tau^2}\varphi_k(\tau)=\zeta_k(\tau),\quad\text{in }(0,R),
\end{equation*}
where
\begin{equation}\label{eq_defzeta}
  \zeta_k(\tau)=\frac{\kappa_s}{\tau^{2-2s}}\int_{{\mathbb
      S}^{N-1}}h(\tau\theta')w(0,\tau\theta') \psi_k(0,\theta')\,dS'-
m^2 \varphi_k(\tau).
\end{equation}
A direct calculation shows that, for some $c_1^k,c_2^k\in\R$,
\begin{equation}\label{eq:varphik}
\varphi_k(\tau)=\tau^{\sigma^+_k}
\bigg(c_1^k+\int_\tau^R\frac{t^{-\sigma^+_k+1}}{\sigma^+_k-\sigma^-_k}
\zeta_k(t)\,dt\bigg)+\tau^{\sigma^-_k}
\bigg(c_2^k+\int_\tau^R\frac{t^{-\sigma^-_k+1}}{\sigma^-_k-\sigma^+_k}
\zeta_k(t)\,dt\bigg),
\end{equation}
with $\sigma^\pm_k=-\frac{N-2s}{2}\pm\sqrt{\big(\frac{N-2s}
    {2}\big)^{2}+\mu_k(a)}$.
From \eqref{eq:ipo1}, \eqref{eq:23}, and  Lemma \ref{l:stima0},
 we
deduce that, for all $i=j_0,\dots,j_0+M-1$,
\begin{equation}\label{eq:zeta}
  \zeta_{i}(\tau)=O(\tau^{-2+\chi+\sigma_{i}^+})\quad\text{as }\tau\to 0^+.
\end{equation}
Consequently, the functions $t\mapsto t^{-\sigma^+_{i}+1}
\zeta_{i}(t)$, $t\mapsto
t^{-\sigma^-_{i}+1}\zeta_{i}(t)$
belong to $L^1(0,R)$. Hence
\[
\tau^{\sigma^+_{i}}
\bigg(c_1^{i}+\int_\tau^R\frac{\rho^{-\sigma^+_{i}+1}}{\sigma^+_{i}-\sigma^-_{i}}
\zeta_{i}(\rho)\,d\rho\bigg)=o(\tau^{\sigma^-_{i}})\quad\text{as }\tau\to0^+,
\]
and then, by \eqref{eq:23}, there must be
\begin{equation*}
c_2^{i}=-\int_0^R\frac{t^{-\sigma^-_{i}+1}}{\sigma^-_{i}-\sigma^+_{i}}
\,\zeta_{i}(t)\,dt.
\end{equation*}
Using (\ref{eq:zeta}), we then deduce that
\begin{align}\label{eq:12}
  \tau^{\sigma^-_{i}}
  \bigg(c_2^{i}+\int_\tau^R\frac{t^{-\sigma^-_{i}+1}}{\sigma^-_{i}-\sigma^+_{i}}
  \zeta_{i}(t)\,dt\bigg)&=\tau^{\sigma^-_{i}}
  \bigg(\int_0^\tau
  \frac{t^{-\sigma^-_{i}+1}}{\sigma^+_{i}-\sigma^-_{i}}
  \zeta_{i}(t)\,dt\bigg)=O(\tau^{\sigma^+_{i}+\chi})
\end{align}
as $\tau\to0^+$. Combining (\ref{eq:varphik}) and
(\ref{eq:12}), we obtain that,  for all $i=j_0,\dots,j_0+M-1$,
\begin{equation}\label{eq:24}
\varphi_i(\tau)=\tau^{\sigma^+_i}
\bigg(c_1^i+\int_\tau^R\frac{t^{-\sigma^+_i+1}}{\sigma^+_i-\sigma^-_i}
\zeta_i(t)\,dt+O(\tau^{\chi})\bigg)\quad\text{as }\tau\to0^+.
\end{equation}
Let us assume by contradiction that
$\lim_{\lambda\to0^+}\lambda^{-2\gamma}H(\lambda)=0$. Then, (\ref{eq:17bis})
would imply that $\lim_{\tau\to0^+}\tau^{-\sigma_{i}^+}\varphi_{i}(\tau)=0$
for all
$i\in\{j_0,\dots,j_0+M-1\}$.
Hence, in view of \eqref{eq:24},
\[
c_1^{i}+\int_0^R\frac{t^{-\sigma^+_{i}+1}}{\sigma^+_{i}-\sigma^-_{i}}
\zeta_{i}(t)\,dt=0,
\]
which, together with (\ref{eq:zeta}), implies
\begin{align}\label{eq:13tris}
  \tau^{\sigma^+_{i}}
  \bigg(c_1^{i}+\int_\tau^R\frac{t^{-\sigma^+_{i}+1}}{\sigma^+_{i}-\sigma^-_{i}}
  \zeta_{i}(t)\,dt\bigg)= \tau^{\sigma^+_{i}}
  \int_0^\tau\frac{t^{-\sigma^+_{i}+1}}{\sigma^-_{i}-\sigma^+_{i}}
  \zeta_{i}(t)\,dt=O(\tau^{\sigma^+_{i}+\chi})
\end{align}
as $\tau\to0^+$. Collecting \eqref{eq:varphik}, (\ref{eq:12}), and
(\ref{eq:13tris}), we conclude that
$\varphi_{i}(\tau)=O(\tau^{\sigma^+_{i}+\chi})$ as $\tau\to0^+$ for every
$i\in\{j_0,\dots,j_0+M-1\}$,
namely,
\[
\sqrt{H(\tau)}\,(w^\tau,\psi)_{L^2({\mathbb
  S}^{N}_+;\theta_1^{1-2s})}=
O(\tau^{\gamma+\chi})\quad\text{as
}\tau\to0^+
\]
for every $\psi\in {\mathcal A}_0=\mathop{\rm span}\{
\psi_i\}_{i=j_0}^{j_0+M-1}$, where ${\mathcal A}_0$ is the eigenspace
of problem \eqref{eq:4} associated to the eigenvalue
  $\mu_{j_0}(a)=\mu_{j_0+1}(a)=\cdots=\mu_{j_0+M-1}(a)$.
  From Lemma \ref{l:stima_N_sopra} part (iv),
there exists $C(\chi)>0$ such that $\sqrt{H(\tau)}\geq
C(\chi)\tau^{\gamma+\frac\chi2}$ for $\tau$ small, and therefore
\begin{equation}\label{eq:26}
(w^\tau,\psi)_{L^2({\mathbb
  S}^{N}_+;\theta_1^{1-2s})}=
O(\tau^{\frac{\chi}2})\quad\text{as
}\tau\to0^+
\end{equation}
for every $\psi\in{\mathcal A}_0$.  From Lemma \ref{l:blowup},
 for every sequence $\tau_n\to0^+$, there exist a subsequence
$\{\tau_{n_k}\}_{k\in\N}$ and an eigenfunction $\widetilde
\psi\in{\mathcal A}_0$
\begin{equation}\label{eq:25}
  \int_{{\mathbb S}^{N}_+}\theta_1^{1-2s}\widetilde\psi^2(\theta)dS=1\quad\text{and} \quad
w^{\tau_{n_k}}\to \widetilde \psi\quad\text{in } L^2({\mathbb
  S}^{N}_+;\theta_1^{1-2s}).
\end{equation}
From (\ref{eq:26}) and (\ref{eq:25}), we infer that
\[
0=\lim_{k\to+\infty}(w^{\tau_{n_k}},\widetilde\psi)_{L^2({\mathbb
  S}^{N}_+;\theta_1^{1-2s}) }
=\|\widetilde\psi\|_{ L^2({\mathbb
  S}^{N}_+;\theta_1^{1-2s})}^2=1,
\]
thus reaching a contradiction.
\QED

\subsection{Proof of Theorem \ref{t:asym}}
Identity
  (\ref{eq:35du}) follows from part (i) of Lemma \ref{l:blowup}, thus
  there exists $k_0\in \N$, $k_0\geq 1$, such that $\gamma=\lim_{r\to
    0^+}{\mathcal N}(r)=-\frac{N-2s}{2}+\sqrt{\big(\frac{N-2s}
    {2}\big)^{2}+\mu_{k_0}(a)}$.  Let us denote as $M$ the
  multiplicity of $\mu_{j_0}(a)$ so that, for some $j_0\in\N$,
  $j_0\geq 1$, $j_0\leq k_0\leq j_0+M-1$,
  $\mu_{j_0}(a)=\mu_{j_0+1}(a)=\cdots=\mu_{j_0+M-1}(a)$
and let $\{
\psi_i\}_{i=j_0}^{j_0+m-1}$ be an $L^2({\mathbb
  S}^{N}_+;\theta_1^{1-2s})$-orthonormal basis for the eigenspace
associated to $\mu_{k_0}(a)$.

Let $\{\tau_n\}_{n\in\N}\subset (0,+\infty)$ such that
$\lim_{n\to+\infty}\tau_n=0$. Then, by Lemma
\ref{l:blowup} part (ii) and Lemmas \ref{l:limite} and \ref{l:limitepositivo},
there exist a subsequence $\{\tau_{n_k}\}_{k\in\N}$ and $M$ real
numbers $\beta_{j_0},\dots,\beta_{j_0+M-1}\in\R$ such that
$(\beta_{j_0},\beta_{j_0+1},\dots,\beta_{j_0+M-1})\neq(0,0,\dots,0)$
and
\begin{align}\label{eq:36}
&\tau_{n_k}^{-\gamma}w(\tau_{n_k}\theta)\to
\sum_{i=j_0}^{j_0+M-1} \beta_i\psi_{i}(\theta)\quad \text{in }
C^{0,\alpha}({\mathbb S}^{N}_+)  \quad \text{as }k\to+\infty,\\
\label{eq:47}&\tau_{n_k}^{-\gamma}w(0,\tau_{n_k}\theta')\to
\sum_{i=j_0}^{j_0+M-1} \beta_i\psi_{i}(0,\theta')\quad \text{in }
C^{1,\alpha}({\mathbb S}^{N-1}) \quad \text{as }k\to+\infty,
\end{align}
for some $\alpha\in(0,1)$.
 We now prove
that the $\beta_i$'s depend neither on the sequence
$\{\t_n\}_{n\in\N}$ nor on its subsequence
$\{\t_{n_k}\}_{k\in\N}$.

   Defining $\varphi_i$ and $\zeta_i$
as in \eqref{eq:23} and \eqref{eq_defzeta}, from
(\ref{eq:36}) it follows that, for any $i=j_0,\dots, j_0+M-1$,
\begin{equation}\label{eq:43}
\tau_{n_k}^{-\gamma}\varphi_i(\tau_{n_k})
\to\beta_i
\end{equation}
as $k\to+\infty$.  As deduced in the proof of Lemma
\ref{l:limitepositivo}, for any $i=j_0,\dots, j_0+M-1$
and $\tau\in(0,R]$ there holds
\begin{align}\label{eq:44}
\varphi_i(\tau)&=
\tau^{\sigma^+_i}
\bigg(c_1^i+\int_\tau^R\frac{t^{-\sigma^+_i+1}}{\sigma^+_i-\sigma^-_i}
\zeta_i(t)\,dt\bigg)+
\tau^{\sigma^-_{i}}
  \bigg(\int_0^\tau
  \frac{t^{-\sigma^-_{i}+1}}{\sigma^+_{i}-\sigma^-_{i}}
  \zeta_{i}(t)\,dt\bigg)\\
\notag &=
\tau^{\sigma^+_i}
\bigg(c_1^i+\int_\tau^R\frac{t^{-\sigma^+_i+1}}{\sigma^+_i-\sigma^-_i}
\zeta_i(t)\,dt+O(\tau^{\chi})\bigg)\quad\text{as }\tau\to0^+,
\end{align}
for some $c_1^i\in\R$.
Choosing $\tau=R$ in the first line of (\ref{eq:44}), we obtain
\[
c_1^i=R^{-\sigma^+_i}\varphi_i(R)-R^{\sigma^-_i-\sigma^+_i}\int_0^R
  \frac{s^{-\sigma^-_i+1}}{\sigma^+_i-\sigma^-_i}
  \zeta_i(s)\,ds.
\]
Hence (\ref{eq:44}) yields
\[
\tau^{-\gamma}\varphi_i(\tau)\to
R^{-\sigma^+_i}\varphi_i(R)-R^{\sigma^-_i-\sigma^+_i}\int_0^R
  \frac{t^{-\sigma^-_i+1}}{\sigma^+_i-\sigma^-_i}
  \zeta_i(t)\,dt+\int_0^R\frac{t^{-\sigma^+_i+1}}{\sigma^+_i-\sigma^-_i}
\zeta_i(t)\,dt
\]
as $\tau\to0^+$,
and therefore from (\ref{eq:43}) we deduce that \eqref{eq:betai}
holds; in particular the $\beta_i$'s depend neither on the sequence
$\{\tau_n\}_{n\in\N}$ nor on its subsequence
$\{\tau_{n_k}\}_{k\in\N}$, thus implying that the convergences in
(\ref{eq:36}) and \eqref{eq:47}
 actually hold as $\tau\to 0^+$
and proving the theorem.\QED

We are now in position to prove Theorem \ref{t:asym-frac} and its
corollaries.

\smallskip\noindent {\bf Proof of Theorem \ref{t:asym-frac}.}
Let $u\in H^s(\R^N)$ be a nontrivial weak solution to $Hu=0$ in
$\Omega$. By Theorems \ref{th:extmag} and \ref{th:ext} in the
appendices there exists a unique
$w=\mathcal H(u)\in H^1(\R^{N+1}_+;t^{1-2s})$ weakly solving
\[
\begin{cases}
  -\dive(t^{1-2s}\nabla w)+m^2 t^{1-2s}w=0,&\text{in }\R^{N+1}_+,\\
w=u,&\text{on }\partial \R^{N+1}_+=\{0\}\times \R^N,
\end{cases}
\]
which also satisfies
$$
-\lim_{t\to 0^+}t^{1-2s}\frac{\partial w}{\partial
  t}(x)=\kappa_s(-\Delta+m^2)^su(x), \quad\text{in }H^{-s}(\R^N)
$$
in a weak sense. Therefore $w$ solves \eqref{eq:wHextended} in the
sense of  \eqref{eq:wH8}. Then Theorem \ref{t:asym-frac} follows from
Theorem \ref{t:asym}.\QED

\smallskip\noindent {\bf Proof of Corollary
  \ref{t:lambda0-asym-frac}.} It follows as a particular case of
Theorem \ref{t:asym-frac}  in the case $a\equiv 0$.\QED

\smallskip\noindent {\bf Proof of Theorem
  \ref{t:sun}.} It follows from Theorem \ref{t:asym-frac}, observing
that if, by contradiction, $u\not\equiv 0$, then convergences stated
Theorem \ref{t:asym-frac} would hold, thus contradicting that
$u(x)=o(|x|^n)$ as $|x|\to 0$ if $n>-\frac{N-2s}{2}+\sqrt{\big(\frac{N-2s}
    {2}\big)^{2}+\mu_{k_0}(a)}$.\QED

\smallskip\noindent {\bf Proof of Theorem
  \ref{t:unspm}.} The proof follows from Corollary
\ref{t:lambda0-asym-frac} arguing as in the proof of \cite[Theorem 1.4]{FF}.
\QED

 \section{Appendix A: Extension theorem}\label{s:ext-th}

Let $s\in (0,1)$ and $N\in \N^*$. Throughout this section
$\RNp:=\{z=(t,x)\,:\,t>0,\, x\in\R^N\}$.
Let $P(D)= P(D_x)$ be a pseudo-differential operator  with constant  coefficients and
Fourier transform (symbol) $P(\xi)\geq0$ with order $\ell\in\R$. We mean $|P(\xi)|\leq C(1+ |\xi|)^\ell$, for some positive constant $C$.  For every $s\in(0,1)$,
define the $s$-power of  $P(D)$ as
$$
\widehat{P(D)^s u}(\xi)= P(\xi)^{s} \widehat{u}(\xi).
$$
Assume  that the bilinear form
$$
(u,v)\mapsto \int_{\R^N}(P(\xi))^{s}\widehat{u}\overline{\widehat{v}}\,d\xi =\int_{\R^N}u(P(D))^sv\,dx
$$
defines a scalar product in $C^\infty_c(\R^N)$ for every $s\in(0,1]$.  Let
 $\dot{H}^s_{D}(\R^N)$  be the  completion of $C^\infty_c(\R^N)$ with respect to the
 above  scalar product.
 Next we
  define the space  $\dot{H}_{D}^1(\RNp;t^{1-2s})$ to be the completion
  of $C^\infty_c(\ov{\RNp})$ with respect to the norm
$$
\bigg( \int_{\RNp}t^{1-2s} w P(D)w\,dxdt+\int_{\RNp}t^{1-2s}\left(\frac{\de w}{ \de t}\right)^{\!\!2}\,dxdt\bigg)^{\!\!1/2}  .
$$
Scalar products in the above spaces are denoted
as $\la \cdot, \cdot\ra_{ \dot{H}^s_{D}(\R^N)}$ and  $\la\cdot,\cdot  \ra_{\dot{H}_{D}^1(\RNp;t^{1-2s} )} $.\\
Under the above setting and  assumptions, the following result holds.
\begin{Theorem}\label{th:extmag}
 Let $s\in(0,1)$ and $u\in \dot{H}^s_{D}(\R^N)$. Then there exists a unique
  $w\in \dot{H}_{D}^1(\RNp;t^{1-2s})$  solution to the problem
\be\label{eq:extmag}
\begin{cases}
 t^{1-2s} P(D) w - (t^{1-2s}w_t )_t =0,&\textrm{ in }\RNp,\\
w=u,&\textrm{ on }\R^N,
\end{cases}
\ee
where the subscript $t$  means derivatives with respect to $t$.
In addition
\be\label{eq:PDs}
 -\lim_{t\to 0} t^{1-2s}\frac{\de w}{\de t}=\k_s\,  (P(D) )^{s}  u\quad \textrm{ in } \dot{H}^{-s}_{D}(\R^N),
\ee in the sense that: for any $\Psi\in \dot{H}_{D}^1(\RNp;t^{1-2s}
) $
$$
\la w,\Psi \ra_{\dot{H}_{D}^1(\RNp;t^{1-2s} )}=\k_s\la u, \Psi\ra_{
\dot{H}^s_{D}(\R^N)}.
$$
Here $\dot{H}^{-s}_{D}(\R^N) $ denotes the dual of
$\dot{H}^s_{D}(\R^N)$ while
\be\label{eq:kappas}
\k_s=2^{1-2s}\frac{\G(1-s)}{\G(s)}
\ee
and $\G$ is the usual Gamma function.
\end{Theorem}
Extension theorems found useful applications in the study of
fractional partial differential equations. For $P(D)=-\D$, see
\cite{CSilv}. We also quote \cite{ST} with $P(D)$ a second order
differential operator with possibly non constant coefficients, see
also \cite{CG}. A main point in our result is that the function space is explicitly given.

\subsection{Proof of Theorem \ref{th:extmag}}\label{s:PT}
We start with some preliminaries. For any  $v\in C^\infty_c(\R^N)$,
we define   $\calH(v)$ via its Fourier transform with respect to the
variable $x$ as $\widehat{
\calH(v)}(t,\xi)=\widehat{v}(\xi)\vartheta(\sqrt{P(\xi)}t)$, where
$\vartheta \in H^1(\R_+; t^{1-2s})$ solves   the ordinary
differential equation: \be\label{eq:Bess}
\begin{cases}
\vartheta''+\frac{(1-2s)}{t}\vartheta'-\vartheta=0,\\
\vartheta(0)=1.
\end{cases}
\ee
 We note that $\vartheta$ is a given by a  Bessel function:
\be
\vartheta(r)=\frac{2}{\G(s)}\left( \frac{r}{2}\right)^s\,K_s(r),
\ee
where, $K_{\nu}$  denotes the modified Bessel function of the second
kind with order $\nu$. It solves the equation
\be\label{eq:BessK}
r^2 K_\nu''+r K_\nu'-(r^2+\nu^2 )K_\nu=0.
\ee
We have, see
\cite{El}, for $\nu>0$,
\be\label{eq:decKnu0}
 K_{\nu}(r)\sim
\frac{\G(\nu)}{2} \left(\frac{r}{2} \right)^{-\nu}
\ee
 as $r\to 0$
and $K_{-\nu}=K_{\nu}$ for $\nu<0$, while
\be\label{eq:decKnuInf}
K_{\nu}(r)\sim \frac{\sqrt{\pi }}{\sqrt{2}} r^{-1/2}e^{-r}\ee as
$r\to +\infty$. By using the identity
 $$
K_\nu'(r)=-\frac{\nu}{r}K_\nu -K_{\nu-1},
$$
we get
\be\label{eq:kappas}
\k_s=\int_0^\infty
t^{1-2s}(|\vartheta'(t)|^2+|\vartheta(t)|^2)\, dt=-\lim_{t\to0}
t^{1-2s}\vartheta'(t)=2^{1-2s}\frac{\G(1-s)}{\G(s)} .
\ee
 Since
$v\in C^\infty_c(\R^N)$, $\widehat{v}$ decays faster
than any polynomial. Then $  \calH(v)\in
\dot{H}_{D}^1(\RNp;t^{1-2s})$ and in addition it satisfies the
equation \be\label{eq:Hvharm}
 t^{1-2s} P(D)\calH(v)  - (t^{1-2s}\calH(v)_t )_t =0\quad\textrm{ in }\RNp.
\ee
We  start by showing that  $P(D)$ satisfies the \textit{trace
property} that any $w\in \dot{H}_{D}^1(\RNp;t^{1-2s}) $ has a trace
which belongs to $\dot{H}^s_{D}(\R^N) $.

\begin{Proposition}\label{prop:trace}
There exists a (unique) linear   trace operator
$$
T: \dot{H}_{D}^1(\RNp;t^{1-2s})\to \dot{H}^s_{D}(\R^N)
$$
 such that
$T(w)(t,x)=w(0,x)$ for any $w\in C^\infty_c(\ov{\RNp})$ and moreover
\be\label{eq:Trace}
 \k_s\, \|T(w)\|_{ \dot{H}^s_{D}(\R^N) }^2\leq \|w\|_{
\dot{H}_{D}^1(\RNp;t^{1-2s}) }^2\quad \textrm{ for all } w\in
\dot{H}_{D}^1(\RNp ;t^{1-2s}), \ee where $\k_s$ is given by
\eqref{eq:kappas}. Equality holds in \eqref{eq:Trace} for some function $w\in\dot{H}_{D}^1(\RNp;t^{1-2s} ) $
if and only if $t^{1-2s} P(D)w  - (t^{1-2s}w_t )_t =0$ in $\RNp$.
\end{Proposition}
\proof
Let $v\in C^\infty_c(\R^N)$. By \eqref{eq:Hvharm},
 we have that   any $w\in C^\infty_c(\ov{\RNp})$ such that $w(0,\cdot)=v$ on $\R^N$ satisfies
\be\label{eq:minenrj}
 \|\calH(v)\|_{\dot{H}_{D}^1(\RNp;t^{1-2s} )}^2\leq  \|w\|_{\dot{H}_{D}^1(\RNp;t^{1-2s} )}^2.
\ee
Thanks to Parseval identity, we have
\begin{align*}
 \|\calH(v)&\|_{\dot{H}_{D}^1(\RNp;t^{1-2s} )}^2=\int_{\RNp}t^{1-2s}P(\xi)  \widehat{\calH(v) }^2d\xi dt +
\int_{\RNp}t^{1-2s}\bigg(\frac{ \de \widehat{ \calH(v)} }{\de
  t}\bigg)^{\!\!2}\, d\xi dt\\
&=\int_{\RNp}t^{1-2s}P(\xi)  \widehat{v}^2(\xi)|\vartheta(\sqrt{P(\xi)}t)|^2 d\xi dt
 +
\int_{\RNp}t^{1-2s}\widehat{v}^2(\xi)\bigg (\frac{\partial}{\partial
  t}\vartheta(\sqrt{P(\xi)}t)\bigg)^{\!\!2}\, d\xi dt\\
&=\int_{\RNp}t^{1-2s}P(\xi) \widehat{v}^2(\xi)|\vartheta(\sqrt{P(\xi)}t)|^2 d\xi dt +
\int_{\RNp}t^{1-2s}\widehat{v}^2 (\xi)P(\xi)
|\vartheta'(\sqrt{P(\xi)}t)|^2  \, d\xi dt\\
&=\bigg(\int_{\R^N} (P(\xi) )^s\widehat{v}^2(\xi) \, d\xi\bigg)\bigg(\int_0^\infty
t^{1-2s}(|\vartheta'(t)|^2+|\vartheta(t)|^2) dt\bigg).
\end{align*}
We conclude that, for any  $v\in C^\infty_c(\R^N)$,
\be\label{eq:Isom}
 \|\calH(v)\|_{\dot{H}_{D}^1(\RNp;t^{1-2s} )}^2=\k_s\|v\|_{ \dot{H}^s_{D}(\R^N)}^2.
\ee
This with  \eqref{eq:minenrj} implies that
$$
\k_s\|w(0,\cdot)\|_{ \dot{H}^s_{D}(\R^N)}^2\leq
\|w\|_{\dot{H}_{D}^1(\RNp;t^{1-2s} )}^2\quad \text{for all } w\in
C^\infty_c(\ov{\RNp}).
$$
The operator $T$ is now defined as the unique extension of the operator $w\mapsto w(0,\cdot)$.

\QED
%%%%%%%%%%%%%%%%%%%%%%%%%%%%%%%%%%%%%%%%%%%%%%%%%%%%%%%%%%%%5

For sake of simplicity, in this paper, we have
denoted the trace of  a function  $w \in \dot{H}_{D}^1(\RNp;t^{1-2s}
)$ with the same letter $w$.

\subsubsection{Proof of Theorem \ref{th:extmag}}
We first consider $u\in C^\infty_c(\R^N)$. In this case $w=\calH(u)$ and it is, of course,
unique in $ \dot{H}_{D}^1(\RNp;t^{1-2s} ) $.

Now we observe that   $-\lim_{t\to 0}t^{1-2s}\frac{\de
\widehat{\calH(u)}(t,\xi)}{\de t} =\kappa_s(P(\xi))^s\widehat{u}(\xi)$ so
that
$$
-\lim_{t\to 0}t^{1-2s}\frac{\de \calH(u)}{\de
t}=\kappa_s(P(D))^su\quad\text{in }  \dot{H}^{-s}_{D}(\R^N).
$$
By \eqref{eq:Hvharm} and
Proposition \ref{prop:trace}, we deduce, after integration by parts,
that
\be\label{eq:sol}
 \la \calH(u),\Psi \ra_{H_{D}^1(\RNp;t^{1-2s}
)}=\k_s\la u, \Psi\ra_{ \dot{H}^s_{D}(\R^N)}
\ee
for any  $\Psi\in  \dot{H}_{D}^1(\RNp;t^{1-2s} )$, and this proves the
theorem in this case.

For the general case $u\in  \dot{H}^s_{D}(\R^N)$, there exists a
sequence $u_n\in C^\infty_c(\R^N)$ such that $u_n\to u$ in
$\dot{H}^s_{D}(\R^N)$. It turns out that $\calH(u_n) \rightharpoonup
\ti{w}$ in $ \dot{H}_{D}^1(\RNp;t^{1-2s} )$ and  $Tr(\calH(u_n))
\rightharpoonup Tr(\ti{w})=u$ in $\dot{H}^s_{D}(\R^N)$ . In
particular for every $\psi\in C^\infty_c(\RNp)$
$$
\la \ti{w},\psi\ra_{ \dot{H}_{D}^1(\RNp;t^{1-2s} )}= 0 .
$$
This implies that $\ti{w}=w$ and it  is unique in $
\dot{H}_{D}^1(\RNp;t^{1-2s} ) $. By \eqref{eq:sol}
$$
\la \calH(u_n),\Psi \ra_{\dot{H}_{D}^1(\RNp;t^{1-2s} )}=\k_s\la u_n,
\Psi\ra_{ \dot{H}^s_{D}(\R^N)}
$$
 for any  $\Psi\in  \dot{H}_{D}^1(\RNp;t^{1-2s} )$.
Taking the limit as $n\to\infty$,  we get the desired result.
\QED

\begin{Remark}
We note that the trace operator  $T$ defined in Proposition
\ref{prop:trace} is surjective. To see that, we argue by density.
Let $v\in  \dot{H}^s_{D}(\R^N)$. There exists a sequence $v_n\in
C^\infty_c(\R^N)$ such that $v_n\to v$ in $  \dot{H}^s_{D}(\R^N)$.
By \eqref{eq:Isom} $\calH(v_n)$ is bounded and thus converges (up to subsequences) weakly
to some function $w\in \dot{H}_{D}^1(\RNp;t^{1-2s} )$ and Theorem
\ref{th:extmag} implies that the convergence is strong   and thus
$$
\| w\|_{\dot{H}_{D}^1(\RNp;t^{1-2s} ) }^2=\k_s \|T(
w)\|_{\dot{H}^s_{D}(\R^N)}^2= \k_s\|v\|_{\dot{H}^s_{D}(\R^N)}^2.
$$
\end{Remark}

\section{Appendix B: The relativistic Schr\"{o}dinger operator $\Dsm$}\label{s:Shro}
 Given $ m>0$, letting $P(D)=-\D+m^2$,  we have $ H^{1}_D(\RNp;t^{1-2s})= H^{1}(\RNp;t^{1-2s}) $
and $H^s_D(\R^N)=H^s(\R^N)$.
Applying  Theorem \ref{th:extmag}, we have the following result.
\begin{Theorem}\label{th:ext}
 Let $u\in H^s(\R^N)$ and let $w\in H^{1}(\RNp;t^{1-2s}) $  be the unique solution  to the problem
\be\label{eq:ext}
\begin{cases}
- \div(t^{1-2s}\n w)+m^2t^{1-2s} w=0,&\text{in }\RNp,\\
w=u, &\text{on }\R^N.
\end{cases}
\ee
Then
$$
 -\lim_{t\to 0} t^{1-2s}\frac{\de w}{\de t}=\k_s\Dsm u\quad \textrm{ in } H^{-s}(\R^N).
$$
\end{Theorem}
\subsection{Bessel Kernel}\label{s:Bess-Kernel}
We can observe that the Bessel kernel $P_m(t,x)$ is given by the
Fourier transform of the mapping $\xi\mapsto
\vartheta(\sqrt{|\xi|^2+m^2}\,t)$, where $\vartheta$ is the Bessel
function solving the differential equation \eqref{eq:Bess}
and yet  we can determine it explicitly.

Let $U$  satisfy
$$
- \div(t^{1-2s}\n U)+m^2t^{1-2s} U=0,\quad\textrm{ in }\RNp.
$$
We have that   $V= t^{1-2s} \frac{\de U}{\de t}$ solves the conjugate
problem:
$$
- \div(t^{-1+2s}\n V)+m^2t^{-1+2s} V=0,\quad \textrm{ in }\RNp.
$$
We look for $F$  (the  fundamental solution) which satisfies
$$
- \div(|t|^{-1+2s}\n F)+m^2|t|^{-1+2s} F=\d_0,\quad \textrm{ in }\R^{N+1}.\\
$$
By direct computations we have
$$
F(z)=C_{N,s}\,
m^{(N+2s-2)/2}|z|^{\frac{-2s-N+2}{2}}K_{\frac{N+2s-2}{2}}(m|z|),
$$
where $C_{N,s} $ is a normalizing constant and  $K_{\nu}$  denotes the modified Bessel function of the second
kind with order $\nu$ solving \eqref{eq:BessK}. Hence the choice of the Bessel Kernel in $\RNp$ is
 \[
P_m(t,x)=-t^{-1+2s}\frac{\de F(t,x)}{\de t}.
\]
Using the identity $K_\nu'(r)=\frac{\nu}{r}K_\nu -K_{\nu+1}$, we
obtain \be\label{eq:Bess-Kernel} {P}_m(z)=C_{N,s}'\,
t^{2s}m^{\frac{N+2s}{2} }
|z|^{-\frac{N+2s}{2}}K_{\frac{N+2s}{2}}(m|z|). \ee
By using \eqref{eq:decKnu0} we deduce  that  $C'_{N,s}$ is given by
$$
C_{N,s}'=p_{N,s}\frac{2^{(N+2s)/2-1}}{\G((N+2s)/2)},
$$
where $p_{N,s} $ is the constant for the (normalized) Poisson Kernel with $m=0$, see \cite{CS}.
We refer to  \cite{BRB}, \cite{CKS} for some Green function estimates
for relativistic killed process.
We also refer to \cite{Ste} for estimates of the Bessel Kernel.

We notice that, since $P_m(t,x)$ is the Fourier transform of
$\xi\mapsto \vartheta(\sqrt{|\xi|^2+m^2}t)$, we have
\be\label{eq:intPm} \int_{\R^N}P_m(t,x)dx=\vartheta(mt). \ee
 Now  we
deduce the norm from the Dirichlet form associated to $\Dsm -
m^{2s}$ via the  Bessel kernel $P_m$. For $s=1/2$ and $N=3$, it
was determined in \cite [Theorem 7.12]{LL}.
\begin{Proposition}\label{Prop:bessK}
For every $u\in H^s(\R^N)$, we have that
$$
\int_{\R^{N}}\left[\left|(-\D+m^2)^{\frac{s}{2}}{u}\right|^2 -
m^{2s}{u}^2 \right]\,dx   =\frac{c_{N,s}}{2}m^{\frac{N+2s}{2} }
\int_{\R^{2N}}\frac{(u(x)-u(y))^2}{|x-y|^{\frac{N+2s}{2}} } K_{
\frac{N+2s}{2}}(m|x-y|)\,dxdy,
$$
with
\be\label{eq:cns}
c_{N,s}=\frac{2^{-(N+2s)/2+1}}{\G((N+2s)/2)} \left(\int_{\R^N}\frac{1-\cos(\xi_1)}{|\xi|^{N+2s}} \right)^{-1}=\frac{2^{-(N+2s)/2+1}}{\G((N+2s)/2)}     \pi^{-\frac
  N2}2^{2s}\frac{\Gamma\big(\frac{N+2s}{2}\big)}{\Gamma(2-s)}s(1-s).
\ee
\end{Proposition}
\proof
We know that
\begin{align*}
 -\int_{\R^N}\lim_{t\to 0} t^{1-2s}\frac{\de \vartheta(\sqrt{|\xi|^2+m^2}t)}{\de t}\widehat{u}^2(\xi) \, d\xi
&=\k_s\int_{\R^{N}}(|\xi|^2+m^2)^s\widehat{u}^2(\xi)\,d\xi\\
&=\k_s \int_{\R^{N}}\left|(-\D+m^2)^{\frac{s}{2}}{u}(x)\right|^2\,dx.
\end{align*}
Or equivalently
\be\label{eq:equivlim}
 -2s\int_{\R^N}\lim_{t\to 0}
\frac{\vartheta(\sqrt{|\xi|^2+m^2}t)-1}{t^{2s}}\widehat{u}^2(\xi) \,
d\xi=\k_s
\int_{\R^{N}}\left|(-\D+m^2)^{\frac{s}{2}}{u}(x)\right|^2\,dx.
\ee
We
now compute the left hand side of the above equality using the
Bessel Kernel $P_m$. Given $t>0$, again by Parseval identity, we
have
\begin{eqnarray*}
\int_{\R^N}
\frac{\vartheta(\sqrt{|\xi|^2+m^2}t)-1}{t^{2s}}\widehat{u}^2(\xi) \, d\xi
&=&\frac{1}{ t^{2s}}\int_{\R^N}\Big(u(x)P_m(t,\cdot)*u(x)-u^2(x)\Big)\,dx,
\end{eqnarray*}
where $ P_m(t,\cdot)*u(x)=\int_{\R^N}u(y)P_m(t,x-y)dy$. We normalize
$P_m$ by putting $\ti{P}_m(t,x)=\frac{1}{\vartheta(tm)}P_m(t,x)$ so
that $\int_{\R^N}\ti{P}_m(t,x)dx=1 $. We therefore have for  $t>0$
\begin{align*}
  \int_{\R^N}&
  \frac{\vartheta(\sqrt{|\xi|^2+m^2}t)-1}{t^{2s}}\widehat{u}^2
  (\xi)\,
  d\xi\\
  &= \frac{\vartheta(tm)}{ t^{2s}} \int_{\R^N}\Big(u(x){\ti{P}_m(t,\cdot)*u(x)-u^2(x)}\Big)\,dx +\frac{1}{ t^{2s}}\int_{\R^N}u^2(x)(\vartheta(tm)-1)\,dx\\
  &=-\frac{\vartheta(tm)}{2 t^{2s}}\int_{\R^{2N}}(u(x)-u(y))^2 \ti{P}_m(t,x-y)\,dy dx+\frac{1}{ t^{2s}}\int_{\R^N}u^2(x)(\vartheta(tm)-1)\,dx\\
  &=-\frac{1}{2 t^{2s}}\int_{\R^{2N}}(u(x)-u(y))^2 {P}_m(t,x-y)\,dy dx
  +\frac{1}{ t^{2s}}\int_{\R^N}u^2(x)(\vartheta(tm)-1)\,dx.
\end{align*}
We conclude that for every $t>0$
\begin{align}
\label{eq:1stlim}
&\int_{\R^N}
\frac{\vartheta(\sqrt{|\xi|^2+m^2}t)-1}{t^{2s}}\widehat{u}^2(\xi) \,
d\xi
=\frac{1}{ t^{2s}}\int_{\R^N}u^2(x)(\vartheta(tm)-1)\,dx\\
 \label{eq:2ntlim}
 &\qquad\qquad\qquad -C_{N,s}m^{\frac{N+2s}{2} }\int_{\R^{2N}}
 \frac{(u(x)-u(y))^2}{(t^2+|x-y|^2)^{\frac{N+2s}{4}} } K_{
   \frac{N+2s}{2}}(m(t^2+|x-y|^2)^{1/2})\,dxdy.
\end{align}
We now have to check that we can pass to the limit as $t\to 0$ under
all the above three integrals. Firstly, we observe that the function
$r\mapsto \frac{\vartheta(r)-1}{r^{2s}}$ is decreasing because $K_s$
is decreasing and thus since $u\in H^s(\R^N) $,
 we deduce from \eqref{eq:equivlim} that
\be\label{eq:to0}
\lim_{t\to 0} \int_{\R^N}
\frac{\vartheta(\sqrt{|\xi|^2+m^2}t)-1}{t^{2s}}\widehat{u}^2(\xi) \,
d\xi=-
\frac{\k_s}{2s}\int_{\R^{N}}\left|(-\D+m^2)^{\frac{s}{2}}{u}(x)\right|^2\,dx.
 \ee
For the same reason, we have that
\be\label{eq:to0m}
\lim_{t\to 0}
\frac{1}{
t^{2s}}\int_{\R^N}u^2(x)(\vartheta(tm)-1)\,dx=-m^{2s}\frac{\k_s}{2s}\int_{\R^N}
u^2(x)\,dx.
\ee
Secondly,  thanks
 to the asymptotics of $K_\nu$, we have that there exist $r,R>0$
such that
$$
|z|^{-\nu}K_{\nu}(m|z|)\leq
\begin{cases}
 C |z|^{-2\nu},& \textrm{for } |x-y|<r,\\
C,&\textrm{for } R\geq |x-y|\geq r,\\
 C |z|^{-2\nu}, &\textrm{for } |x-y|> R,
\end{cases}
$$
where $C$ is a positive constant depending only on $N,s,r,R$ and
$m$. Since $u\in H^s(\R^N)$, we can pass the limit as $t\to 0$ under
the integral in \eqref{eq:2ntlim}. This with \eqref{eq:to0} and
\eqref{eq:to0m} in \eqref{eq:1stlim} yields the result. Finally, to
prove \eqref{eq:cns} we use the precise estimate \eqref{eq:decKnu0}
and comparing with the Dirichlet form in the case $m=0$, see
\cite{DPV}.

\begin{Remark}\label{rem:Dsm}
We first remark from the above result that for every $u\in C^2_c(\R^N)$
$$
\Dsm u(x)={c_{N,s}}m^{\frac{N+2s}{2} }
P.V.\int_{\R^{N}}\frac{u(x)-u(y)}{|x-y|^{\frac{N+2s}{2}} } K_{ \frac{N+2s}{2}}(m|x-y|)\,dy+m ^{2s} u(x).
$$
We observe, using    similar arguments as in the the proof of Proposition \ref{Prop:bessK}, that, for
 $$
 w(t,x)=P_m(t,\cdot)*u,
 $$
 with $u\in C^2_c(\R^N)$, we have that
 \begin{equation*}
\begin{cases}
- \div(t^{1-2s}\n w)(t,x)+m^2t^{1-2s} w(t,x)=0,&\textrm{for all  } (t,x)\in\RNp,\\
 -\lim_{t\to 0} t^{1-2s}\frac{\de w}{\de t}(t,x)
 = \k_s\Dsm u (x),&\textrm{for all } x\in \R^N.
\end{cases}
\end{equation*}
\end{Remark}

\noindent
We now prove the following result.
\begin{Proposition}\label{eq:Poiss-cont}
Let $u\in C(\O)$ such that $\int_{\R^N}(1+|x|^{N+2s})^{-1}|u(x)|dx<+\infty$. Let
 $$
 w(t,x)=\big(P_m(t,\cdot)*u\big)(x).
 $$
Then $$\lim_{t\to 0^+}w(t,x)=u(x)$$  for every  $x\in\O$.
\end{Proposition}
\proof
 We recall from \eqref{eq:intPm} that $\int_{\R^N}
 {P_m}(t,x)dx=\vartheta(tm)$ for all $t>0$.
Let $x_0\in\O$; by continuity, for every $\e>0$, there exist
$t_\e,r_\e>0$ such that $|u(y)-\vartheta(tm)u(x_0)|<\e$
for every $y\in B'_{r_\e}(x_0)$ and $0<t<t_\e$. Then
$$
\begin{aligned}
\bigg| w(t,{x_0})&-\vartheta(tm)u(x_0) \bigg|
=\left|\int_{\R^N}(u(y)-  u(x_0)) {P_m}(t,x_0-y)dy \right| \\
& \leq\int_{|y-x_0|<r_\e}|u(y)- u(x_0) | {P_m}(t,x_0-y)dy+
\int_{|y-x_0|\geq r_\e}|u(y)- u(x_0) | {P_m}(t,x_0-y)dy\\
& \leq\e\vartheta(tm) +  \int_{|y-x_0|\geq r_\e}|u(y)-
u(x_0) | {P_m}(t,x_0-y)dy.
\end{aligned}
$$
Using the fact that
$$
|z|^{-\nu}K_{\nu}(m|z|)\leq C_{\nu}|z|^{-2\nu},
$$
we have that
$$
\begin{aligned}
  \int_{|y-x_0|\geq r_\e} {P_m}(t,x_0-y)dy&\leq C t^{2s}\int_{|y-x_0|\geq  r_\e}\frac{1}{(t^2+|y-x_0|^2)^{\frac{N+2s}2}}dy\\
  &\leq C t^{2s}\int_{|y-x_0| \geq
    r_\e}\frac{1}{(|y-x_0|)^{{N+2s}}}dy.
\end{aligned}
$$
Hence we have
\begin{eqnarray*}
\left| w(t,{x_0})- \vartheta(tm)u(x_0) \right|\leq
\e\vartheta(tm)+ t^{2s}C_{r_\e}\bigg[
\int_{\R^N}\frac{|u(x)|}{1+|x|^{N+2s}} dx+ |u(x_0)| \bigg].
\end{eqnarray*}
Sending  $t\to0^+$ and $\e\to0$ respectively, we get the result.

\QED

\subsection{Harnack and local Schauder estimates for the relativistic Schr\"{o}dinger equation}
Let $f\in L^1_{loc}(B_1')$. We recall that
a solution (resp. subsolution, supersolution)  $u\in H^s(\R^N)$ to the equation
 \be\label{locsolrS}
   \Dsm u=(\textrm{resp. }\leq, \geq )\,\,f\quad \textrm{ in } B_1'
 \ee
 satisfies
 $$
 \int_{\R^N}(|\xi|^2+m^2)^s\widehat{u}\ov{\widehat{\vp}}d\xi=(\textrm{resp. }\leq,
 \geq )\,\,\int_{B_1'}f\vp dx
\quad\text{for all }\vp\in H^s(\R^N),
 $$
 or, equivalently, thanks to Proposition \ref{Prop:bessK},
 $$
 \begin{aligned}
\frac{ c_{N,s}}{2}m^{\frac{N+2s}{2} }
\int_{\R^{2N}}\frac{(u(x)-u(y))(\vp(x)-\vp(y))}{|x-y|^{\frac{N+2s}{2}} } &K_{ \frac{N+2s}{2}}(m|x-y|)\,dxdy\\
+m^{2s}&\int_{\R^N}u\vp dx=(\textrm{resp. }\leq, \geq )\,\,\int_{B_1'}f\vp dx.
 \end{aligned}
 $$
 The following regularity result holds.
\begin{Proposition}\label{prop:HSrS}
 Let $a,b\in L^p(B'_1)$, for some $p>\frac{N}{2s}$.
  \begin{enumerate}
  \item  If $u\in H^s(\R^N)$  satisfies $ \Dsm u \leq a(x) u+b(x)$ in $B_1'$ then $u\in L^\infty_{loc}(B_1')$.
    \item  If $u\in H^s(\R^N)$ is  nonnegative and satisfies $ \Dsm u \geq  a(x) u+b(x)$ in $B_1'$, then
   $$
   \inf_{B'_{1/2}}u+ \|b\|_{L^p(B_1')}\geq C \sup_{B'_{1/2}}u.
   $$
    \item  If $u\in H^s(\R^N)$    satisfies $ \Dsm u =  a(x) u+b(x)$ in $B_1'$,
    then $u\in C^{0,\a}_{loc}(B_1')$ and
   $$
  \|u \|_{ C^{0,\a}(\ov{B_{1/2}'}) }\leq C (  \|u \|_{ L^\infty({B_{3/2}'}) }+
  \|b\|_{L^p(B_1')}),
   $$
 \end{enumerate}
where $C>0$ depends only on $N,s,m, \|a\|_{L^{p}(B_1')}$.
\end{Proposition}
\proof
By Theorem \ref{th:extmag}, if  $u\in H^s(\R^N)$    satisfies
$$
\Dsm u = (\textrm{resp. }\leq, \geq )\,\,  a(x) u+b(x)
\quad\textrm{in } B_1'
$$
then there exist a unique $w\in H^1(\RNp;t^{1-2s})$ such that
$$
\begin{cases}
- \div(t^{1-2s}\n w)+m^2t^{1-2s} w=0,&\textrm{in }\RNp,\\
w=u,&\textrm{on }\R^N,\\
 -\frac{1}{\k_s}\lim_{t\to 0} t^{1-2s}\frac{\de w}{\de t}
 =(\textrm{resp. }\leq, \geq )\,\,   a(x) w+b(x),&\textrm{on } B_1',
 \end{cases}
$$
weakly.
The result then follows from Propositions \ref{lem3.1}, \ref{lem3.2} and  \ref{p:Hold}.
\QED
%%%%%%%%%%%%%%%%%%%%%%%%%%%%%%%%%%%%%%%

\label{References}

\end{document}